\documentclass{amsart}

\usepackage{amssymb}
\usepackage{tikz}
\usepackage{verbatim}
\usepackage{hyperref}
\usepackage[english]{babel}
\usepackage{wasysym}
\usepackage[shortlabels]{enumitem}

\usepackage[colorinlistoftodos]{todonotes}
\usepackage{color, soul}
\usepackage{mathtools}

\newtheorem{theorem}{Theorem}[section]
\newtheorem{lemma}[theorem]{Lemma}
\newtheorem{corollary}[theorem]{Corollary}

\theoremstyle{definition}
\newtheorem{definition}[theorem]{Definition}
\newtheorem{notation}[theorem]{Notation}
\newtheorem{example}[theorem]{Example}

\newtheorem{proposition}[theorem]{Proposition}
\newtheorem{fact}[theorem]{Fact}
\newtheorem{observation}[theorem]{Observation}
\newtheorem{claim}[theorem]{Claim}
\newtheorem{conjecture}[theorem]{Conjecture}

\theoremstyle{remark}
\newtheorem{remark}[theorem]{Remark}

\numberwithin{equation}{section}

\newcommand{\aut}{\operatorname{Aut}}

\newcommand{\sym}{\operatorname{Sym}}
\newcommand{\id}{\operatorname{id}}

\newcommand{\csp}{\operatorname{CSP}}
\newcommand{\dom}{\operatorname{Dom}}
\newcommand{\rng}{\operatorname{Rng}}
\newcommand{\rk}{\operatorname{rk}}

\newcommand{\ar}{\operatorname{ar}}
\def\acts{\curvearrowright}

\newcommand{\acl}{\operatorname{acl}}

\newcommand{\tp}{\operatorname{tp}}
\newcommand{\fa}{\mathfrak{A}}
\newcommand{\fb}{\mathfrak{B}}
\newcommand{\fc}{\mathfrak{C}}
\newcommand{\fd}{\mathfrak{D}}

\newcommand{\sign}{\Sigma}

\newcommand{\ms}{\mathcal{M}}

\newcommand{\bigeq}{\nabla}
\newcommand{\smalleq}{\Delta}
\newcommand{\partition}{P}

\newcommand{\gone}{\mathcal{H}}

\newcommand{\perm}{\mathcal{P}}

\newcommand{\smallgroup}{\mathcal{N}}
\newcommand{\biggroup}{\mathcal{S}}
\newcommand{\group}{\mathcal{G}}

\newcommand{\fhom}{\mathcal{FH}}
\newcommand{\fbh}{\mathcal{FBH}}

\newcommand{\typestruct}{\Delta}
\newcommand{\same}{\approx}
\newcommand{\notsame}{\not\approx}
\newcommand{\age}{\operatorname{Age}}
\newcommand{\forb}{\operatorname{Forb}}
\newcommand{\Min}{\operatorname{Min}}
\newcommand{\single}{\mathcal{E}}
\newcommand{\width}{w}

\newcommand{\orb}{o}
\newcommand{\oi}{o^i}
\newcommand{\os}{o^s}

\newcommand\comp[1]{{#1}^{c}}

\begin{document}

\title{Classification of $\omega$-categorical monadically stable structures}
\author{Bertalan Bodor\\ \today}
\thanks{The author has received funding from the European Research Council (Grant Agreement no. 681988, CSP-Infinity).}
\maketitle

\begin{abstract}
	A first-order structure $\fa$ is called \emph{monadically stable} iff every expansion of $\fa$ by unary predicates is stable. In this article we give a classification of the class $\ms$ of $\omega$-categorical monadically stable structures in terms of their automorphism groups. We prove in turn that $\ms$ is smallest class of structures which contains the one-element pure set, closed under isomorphisms, and closed under taking finitely disjoint unions, infinite copies, and finite index first-order reducts. Using our classification we show that every structure in $\ms$ is first-order interdefinable with a finitely bounded homogeneous structure. We also prove that every structure in $\ms$ has finitely many reducts up to interdefinability, thereby confirming Thomas' conjecture for the class $\ms$.
\end{abstract}

\section{Introduction}

	A first-order structure $\fa$ is called \emph{monadically stable} iff every expansion of $\fa$ by unary predicates is stable. We know that this property is also equivalent to a property called  \emph{tree-decomposability} introduced in~\cite{baldwin1985second}. The equivalence follows from the results of~\cite{baldwin1985second} and they are discussed in the introduction in~\cite{lachlan1992}. Further equivalent conditions for monadic stability are given in~\cite{lachlan1992}.
	
	In this article we are interested in those monadically stable structures which are $\omega$-categorical and countable. We denote the class of these structures by $\ms$. Lachlan~\cite{lachlan1992} gave a bottom-up description of the class $\ms$ in terms of automorphism groups. Namely he showed that a structure is contained in $\ms$ if and only if its automorphism group can be obtained from finite structures by a finite number of iterations of certain combinatorial constructions, called \emph{$\omega$-stretches} and \emph{loose unions}. Another interesting characterization of $\ms$ is given by Braunfeld~\cite{braunfeld2019monadic} in terms of orbit growth on subsets: a countable $\omega$-categorical structure is monadically stable if and only if it is stable and the its orbit growth on subsets is slower than exponential (see Theorem~\ref{sam} for the precise formulation).

	In the present article we refine the description of automorphism groups of structures in $\ms$ given in~\cite{lachlan1992} (Section~\ref{sect:class}). We define recursively a sequence of classes of permutation groups $\gone_n'\colon n\in \omega$ such that $\gone_n'$ is exactly the class of automorphism groups of those structures in $\ms$ whose Morley rank are at most $n$. We know from the results of~\cite{lachlan1992} that every monadically stable structure has finite Morley rank. This implies the following.

\begin{theorem}\label{main:str}
	A countable $\omega$-categorical structure $\fa$ is monadically stable if and only if $\aut(\fa)\in \gone':=\bigcup_{n=0}^\infty\gone_n'$.
\end{theorem} 

	The main difference between our description and the one given in~\cite{lachlan1992} is that in our description the groups in $\gone_n'$ have an explicit definition using the groups in $\gone_{n-1}'$.

	We also show that $\gone'$ is the smallest class of permutation groups which contains all finite domain permutation groups, is closed under isomorphisms, and closed under taking finite direct products, wreath products with $\sym(\omega)$, and finite index supergroups. This result also gives rise to a rather simple description of the class $\ms$ without any mention of automorphism groups.

\begin{theorem}\label{desc_struct_intro}
	$\ms$ is the smallest class of structures which contains the one-element pure set, is closed under isomorphisms, and closed taking finitely disjoint unions, infinite copies, and finite index first-order reducts.
\end{theorem}

\subsection*{Finite homogenizability and finite boundedness}

	A structure $\fa$ is called \emph{homogeneous} if every homomorphism between finitely generated substructures of $\fa$ can be extended to an automorphism of $\fa$. A structure is called \emph{finitely homogenizable} if it is first-order interdefinable with a homogeneous structures with a finite relational signature. In~\cite{lachlan1992} Lachlan showed that all structures in $\ms$ are finitely homogenizable.
	
	In the study of finitely homogeneous structures we are particularly interested in the ones which are \emph{finitely bounded} i.e., the ones which can be described by finitely many forbidden structures. The precise definition can be found in Section~\ref{sect:hom}. In Section~\ref{sect:fin_bounded} we show that every structure in $\ms$ is first-order interdefinable with a finitely bounded homogeneous structure (Theorem~\ref{thm:ms_fbh}). This strengthens the aforementioned result by Lachlan.
\subsection*{Connection to CSPs}

	Let $\fb$ be a structure with finite relational signature. Then the \emph{Constraint Satisfaction Problem} (\emph{CSP}) over $\fb$, denoted by $\csp(\fb)$, is the computational problem of deciding whether a given finite structure $\fa$ with the same signature as $\fb$ has a homomorphism to $\fb$. Using concepts and techniques from universal algebra, Bulatov and Zhuk proved that for finite structures $\fb$ the computational complexity of $\csp(\fb)$ satisfies a dichotomy: it is either in $\mathbf{P}$ or it is $\mathbf{NP}$-complete (\cite{BulatovFVConjecture},\cite{ZhukFVConjecture}). The universal algebraic approach for CSPs can also be generalized for countable $\omega$-categorical structures. However as opposed to the finite domain case one cannot expect such a CSP dichotomy for $\omega$-categorical structures since in general we do not even know whether $\csp(\fb)$ is in $\mathbf{NP}$. In fact we know that $\csp(\fb)$ may be undecidable even for homogeneous finitely related structures (see for instance~\cite{BodirskyNesetrilJLC}). One way to guarantee that $\csp(\fb)$ is at least in $\mathbf{NP}$ is to assume that $\fb$ is homogeneous and finitely bounded, or a reduct of such a structure. This motivates the following complexity dichotomy conjecture, originally presented in~\cite{BPP-projective-homomorphisms}.

\begin{conjecture}[Infinite domain CSP dichotomy conjecture]\label{conj:csp}
	Let $\fa$ be a countable finitely bounded homogeneous structure, and let $\fb$ be a finitely related first-order reduct of $\fa$. Then the $\csp(\fb)$ is either in $\mathbf{P}$ or it is $\mathbf{NP}$-complete.
\end{conjecture}

	The conjecture above has been solved for many classes of structures, such as the reducts of
\begin{itemize}
\item the pure set $(\mathbb{N};=)$~\cite{ecsps},
\item the universal linear order $(\mathbb{Q};<)$~\cite{tcsps-journal},
\item the random graph~\cite{BodPin-Schaefer-both},
\item the random poset~\cite{posetCSP18},
\item unary $\omega$-categorical structures~\cite{BodMot-Unary},
\end{itemize}

	and all first-order expansions of the homogeneous RCC5-structure~\cite{bodirsky2020canonical}.

	Theorem~\ref{thm:ms_fbh} implies that all structures in $\ms$ fall into the scope of the conjecture above. We mention one more property of the class $\ms$ which makes it a good candidate for a CSP dichotomy result, namely that it is closed under taking model-complete cores. A countable $\omega$-categorical structure is called a \emph{model-complete core} if the closure of its automorphism group is the same as its endomorphism monoid. We know that every $\omega$-categorical structure $\fb$ is homomorphically equivalent to an (up to isomorphism unique) model-complete core $\fc$~\cite{Cores-journal}. Moreover the structure $\fc$ is also $\omega$-categorical, and it is called \emph{the model-complete core} of $\fb$. Clearly every countable $\omega$-categorical has the same CSP as its model-complete core. Moreover some of the universal algebraic techniques only work under the assumption that the domain structure is a model-complete core. Thus in the analysis of CSPs we usually prefer to work with the model-complete cores of structures. This reduction is a key step in the solution of all the dichotomy results mentioned above. For this reason when formulating a CSP dichotomy for a class $\mathcal{C}$ of $\omega$-categorical structures it is useful to have a class $\mathcal{C}$ which is closed under taking model-complete cores. This property is satisfied by $\ms$ for simple orbit growth reasons (see Section~\ref{sect:class}).
	
	For more details on infinite domain CSPs we refer the reader to~\cite{Topo},\cite{Bodirsky-HDR} or~\cite{Book}.

\subsection*{Thomas' conjecture}

	In~\cite{RandomReducts} Simon Thomas made the conjecture that a countable homogeneous structure over a finite relational language has only finitely many reducts up to first-order interdefinability. The conjecture has been verified for many well-known homogeneous structures. The list includes the rationals with the usual ordering~\cite{Cameron5}, the countably infinite random graph~\cite{RandomReducts}, the homogeneous universal $K_n$-free graphs~\cite{Thomas96}, the expansion of $({\mathbb Q};<)$ by a constant~\cite{JunkerZiegler}, the universal homogeneous partial order~\cite{Poset-Reducts}, the random ordered graph~\cite{42}, and many more~\cite{agarwal,AgarwalKompatscher,BodJonsPham,BBPP18}. Note that if we drop the assumption that the signature of the homogeneous structure $\fa$ is relational, then the conjecture of Thomas is false even if we keep the assumption that $\fa$ is $\omega$-categorical: already the countable atomless Boolean algebra has infinitely many first-order reducts~\cite{BodorCameronSzabo}.

	Even though Thomas' conjecture has been verified for several individual structures, there does not seem to be much progress towards proving it in its full generality. In a recent paper~\cite{bodirsky2018structures} the authors showed that Thomas' conjecture holds for all cellular structures. In this article we further generalize this result by showing that Thomas' conjecture also holds for the entire class $\ms$ (Theorem~\ref{thomas}). Finally, we mention a recent result by Simon~\cite{simon2018nip} where he showed that Thomas' conjecture holds for another robust class of structures, namely the class of \emph{$\omega$-categorical primitive NIP structures with thorn rank 1}. This class contains for instance all the \emph{generalized random permutations} i.e., the Fra\"{i}ss\'{e} limit $(M;<_1,\dots,<_n)$ of all finite structures equipped with $n$ linear orders.
	
\subsection*{Acknowledgments}

	I would like to thank Andr\'{e}s Aranda for a careful reading and useful feedback on earlier versions of this article, and Manuel Bodirsky for discussion and helpful comments.

\section{Preliminaries}

\subsection{Permutation groups}

	We say that a group $G$ is a \emph{permutation group} if $G$ a subgroup of $\sym(X)$ for some $X$. Note that whenever $G$ is a permutation group then the set $X$ witnessing this fact is unique. We call this set the \emph{domain} of $G$, we denote it by $\dom(G)$.

	If $G\subset\sym(X)$ is a permutation group, and $e\colon X\to Y$ is any injective function, then the map defined as $\iota(e)\colon\gamma\mapsto e\gamma e^{-1}$ is a homomorphism from $G$ to $\sym(\rng(e))$, and is referred to as the homomorphism \emph{induced by $e$}. 

	Two permutation groups $G$ and $H$ are called \emph{isomorphic as permutation groups} if there is a bijection $e\colon\dom(G)\rightarrow \dom(H)$ such that the homomorphism $\iota(e)$ induced by $e$ is an isomorphism between $G$ and $H$. Note that this condition is strictly stronger than saying that $G$ and $H$ are isomorphic \emph{as groups}. For instance the trivial permutation groups $\id(X)$ are isomorphic for every set $X$ as groups, but $\id(X)$ and $\id(Y)$ are only isomorphic as permutation groups if $|X|=|Y|$. All the groups in this article will be permutation groups, and unless otherwise specified, and when we say that two permutation groups are isomorphic we always mean that they are isomorphic \emph{as permutation groups}.

	Assume that $G$ is a permutation group, and $Y\subset \dom(G)$. Then we use the following notation.
\begin{itemize}[noitemsep,topsep=0pt]
\item $G_{Y}$ denotes the \emph{pointwise stabilizer} of $Y$, that is, $G_Y=\{g\in G\colon\forall y\in Y(g(y)=y)\}$.
\item $G_{\{Y\}}$ denotes the \emph{setwise stabilizer} of $Y$, that is, $G_{\{Y\}}=\{g\in G\colon\forall y\in Y(g(y)\in Y)\}$.
\item $G|_Y$ denotes the \emph{restriction of $G$ to $Y$}, that is, $G|_Y=\{h|_Y\colon h\in G\}$, provided that $Y$ is invariant under $G$.
\item $G_{(Y)}:=G_{\{Y\}}|_Y=G_{\{\dom(G)\setminus Y\}}|_Y$.
\item $G_{((Y))}:=G_{\dom(G)\setminus Y}|_Y$.
\end{itemize}

	Whenever the set $Y$ is finite, say $Y=\{y_1,\dots,y_n\}$, we simply write $G_{y_1,\dots,y_n}$ for $G_Y$.

	If $E$ is an equivalence relation defined on a set $X$ then we denote by $X/E$ the set of equivalence classes of $E$, and for $x\in X$ we denote by $[x]_E$ the unique equivalence class of $E$ containing $x$. Then by definition $(x,y)\in E$ iff $[x]_E=[y]_E$. An equivalence relation $E$ defined on $\dom(G)$ for some permutation group $G$ is called a \emph{congruence} of $G$ if and only if $E$ is preserved by every permutation in $G$. In this case there is a natural map $G\to \sym(\dom(G)/E)$ that maps $g\mapsto g/E$, where $g/E([x]_E)=[g(x)]_E$. We define $G/E$ as $\{g/E\colon g\in G\}$.

\subsection{Topology}

	Every permutation group is naturally equipped with a topology, called the \emph{topology of pointwise convergence}. This topology can be defined as the subspace topology of the product space ${}^XX$ where $X$ is equipped with the discrete topology. A sequence $(g_n)_n$ in $G$ converges to a permutation $g\in \sym(X)$ with respect to this topology if and only if for all finite subsets $F$ of $X$ there exists an $n$ such that $g_i|_F=g|_F$ for all $i\geq n$.

	We say that a permutation group $G\leq \sym(X)$ is \emph{closed} if it is closed in the topology of pointwise convergence. Equivalently, $G\leq\sym(X)$ is closed if and only if it satisfies the following property: for all $g\in\sym(X)$, if for every finite $F\subset X$ there exists $g'\in G$ such that $g'|F=g|F$ then $g\in G$.

\subsection{Orbit growth functions, oligomorphic groups}

	There are three natural sequences counting orbits attached to a permutation group introduced and discussed in general in~\cite{CameronCounting,Oligo}.

\begin{definition}\label{def:growth}
	Let $G\subseteq \sym(X)$ be a permutation group, and let $n\in \omega$. Then 
\begin{itemize}
\item $\orb_n(G)$ denotes the \emph{number of $n$-orbits of $G$} i.e., the number of orbits of the natural action $G\acts X^n$,
\item $\oi_n(G)$ denotes the \emph{number of injective $n$-orbits of $G$} i.e., the number of those $n$-orbits of $G$ which contain tuples with pairwise different entries,
\item $\os_n(G)$ denotes the \emph{number of orbits of $n$-subsets of $G$} i.e., the number of orbits of the natural action $G\acts {X\choose n}=\{Y\subset X\colon|Y|=n \}$.
\end{itemize}
	A permutation group is called \emph{oligomorphic} if and only if $\orb_n(G)$ is finite for all $n$. 
\end{definition}

\begin{remark}
	It follows easily from the definitions that for any permutation group $G$ the following inequalities hold: $\os_n(G)\leq \oi_n(G)\leq \orb_n(G)\leq n!\os_n(G)$. Therefore the following are equivalent.
\begin{itemize}
\item $G$ is oligomorphic.
\item $\oi_n(G)$ is finite for all $n$.
\item $\os_n(G)$ is finite for all $n$.
\end{itemize}
\end{remark}

	We often use the following easy observation.

\begin{observation}\label{fewer_orbits}
	If $G\leq H$ are permutation groups, then $f_n(H)\leq f_n(G)$ where $f$ is any of the operators $\orb,\oi$ or $\os$. In particular every group containing an oligomorphic group is oligomorphic. 
\end{observation}

\subsection{Direct products}\label{sect:direct}
	
	Let $\{X_i\colon i\in I\}$ be a set of pairwise disjoint sets, and for each $i\in I$ let $G_i$ be a permutation group acting on $X_i$. Then we define the \emph{direct product} of the permutation groups $G_i$, denoted by $\prod_{i\in I}G_i$, to be the set of all permutations of $X:=\bigcup\{X_i\colon i\in I\}$ which can be written as $\bigcup\{f(i)\colon i\in I\}$, for some $f\colon I\to\prod_{i\in I}G_i$ that satisfies $f(i)\in G_i$ for all $i\in I$. Then $\prod_{i\in I}G_i$ as a group is also the usual direct product of the groups $G_i\colon i\in I$.
	
	If $G\leq \sym(X)$ is a permutation group, and $S$ is a set then we define the power $G^S$ to be $\prod_{s\in S}G_s$, where $G_s:=\{(g,s)\colon g\in G\}$ and $(g,s)(g',s)=(gg',s)$. This makes sense since the sets $\dom(G_i)=X\times \{i\}$ are pairwise disjoint. Then $\dom(G^S)=X\times S$, and since all $G_i$ are isomorphic to $G$ it follows that $G^S$ is isomorphic to the power $G^{|S|}$ as a group.

	It follows easily from the definition that any direct product of closed groups is again closed. This also means that any power of a closed group is closed.

	It is easy to see that every $n$-orbit of a direct product $G:=\prod_{i\in I}G_i$ can be written as the union of $n_i$-orbits of $G_i$ for some sequence of natural numbers $(n_i)_{i\in I}$ with $\sum_{i\in I}n_i=n$. This implies that if $I$ is finite then the number of $n$-orbits of $G$ is at most ${|I| \choose n}M$ where $M$ is the maximum of $\orb_n(G_i)\colon i\in I$. In particular if $\orb_n(G_i)$ is finite for all $i\in I$ and $I$ is finite then $\orb_n(G)$ is also finite. Thus any finite direct product of oligomorphic groups is again oligomorphic. 

\subsection{Wreath products}\label{sect:wreath}
	
	Let $G\leq \sym(X)$ and $H\leq \sym(Y)$ be permutation groups. Then we define the \emph{wreath product} of the groups $G$ and $H$, denoted by $G\wr H$, as follows. The domain set of $G\wr H$ is $X\times Y$ and a permutation $\sigma\in \sym(X\times Y)$ is contained in $G\wr H$ if and only if there exist $\beta=\pi(\sigma)\in H$ and $\alpha_y\in G$ such that $\sigma(x,y)=(\alpha_y(x),\beta(y))$. Note that in this case the map $\sigma\mapsto \pi(\sigma)$ defines a surjective homomorphism from $G\wr H$ to $H$, and the kernel of this homomorphism is exactly $G^Y$. Moreover the map $e\colon \beta\mapsto ((x,y)\mapsto (x,\beta(y))$ splits the homomorphism $\pi$, that is $\pi\circ e=\id(Y)$. This implies that $G\wr H$ can be written as a semidirect product $G^Y\rtimes H$ (or more precisely $G^Y\rtimes e(H)$).

	In this document we are only interested in the case when the group $H$ above is the full symmetric group.

\begin{lemma}\label{wr_closed}
	Let $G\leq \sym(X)$ be a permutation group, and $Y$ be a set. Then
\begin{enumerate}
\item if $G$ is closed then so is $G\wr \sym(Y)$, and
\item if $G$ is oligomorphic then so is $G\wr \sym(Y)$.
\end{enumerate}
\end{lemma}

\begin{proof}
	(1) A permutation $\sigma$ of $X\times Y$ is contained in $\sym(X)\wr \sym(Y)$ if and only if it preserves the equivalence relation $E:=\{(x_1,y_1),(x_2,y_2)\colon y_1=y_2\}$. This implies immediately that the group $\sym(X)\wr \sym(Y)$ is closed. Now let $(g_n)_n$ be a sequence in $G\wr \sym(Y)$ converging to a permutation $\sigma\in \sym(X\times Y)$. Then there exist $\beta_n\in \sym(Y),\alpha_{y,n}\in G$ such that $g_n(x,y)=(\alpha_{y,n}(x),\beta_n(y))$. Since $\sym(X)\wr \sym(Y)$ is closed we know that $g\in \sym(X)\wr \sym(Y)$, and thus there exist $\beta\in \sym(Y)$ and $\alpha_y\in \sym(X)$ such that $g(x,y)=(\alpha_y(x),\beta(y))$. We claim that the sequence $(\alpha_{n,y})_n$ converges to $\alpha_y$ for all $y\in Y$. Indeed let $F$ be an arbitrary finite subset of $X$. Then since $(g_n)_n\rightarrow g$ it follows that we can choose an index $n\in \omega$ such that for all $i\geq n$ we have $g|_{F\times \{y\}}=(g_i)|_{F\times \{y\}}$, and thus $\alpha_{i,y}|_F=\alpha_y|_F$. We have shown that $(\alpha_{n,y})_n$ converges to $\alpha_y$, and thus $\alpha_y\in G$ since $G$ is closed. This shows that $g\in G\wr \sym(Y)$.

	(2) Let $n\in \omega$, and let $Y'$ be any $n$-element subset of $Y$. Then every $n$-orbit of $G\wr \sym(Y)$ contains a tuple in $X\times Y'$, and two $n$-tuples from $X\times Y'$ are contained in the same orbit of $G\wr \sym(Y')$ if and only if they are contained in the same orbit of $G\wr \sym(Y)$. Therefore $\orb_n(G\wr \sym(Y))=\orb_n(G\wr \sym(Y'))$. We have already seen that the group $G^{Y'}$ is oligomorphic. On the other hand $G\wr \sym(Y')$ contains $G^{Y'}$, and thus by Observation~\ref{fewer_orbits} it follows that $G\wr \sym(Y')$ is oligomorphic. In particular $\orb_n(G\wr \sym(Y'))=\orb_n(G\wr \sym(Y))$ is finite. Since this holds for all $n\in \omega$ it follows that $G\wr \sym(Y)$ is oligomorphic.
\end{proof}
	
\subsection{Reducts of $\omega$-categorical structures}

	In this article every structure is assumed to be countable.

\begin{notation}
	For a structure $\fa$ we denote by $\sign(\fa)$ the signature of $\fa$.
\end{notation}

	We say that a structure is \emph{relational} if its signature does not contain any function symbols.

\begin{notation}
	Let $\fa$ be a relational structure. Then we denote by $m(\fa)$ the maximal arity of relations of $\fa$. $m(\fa)$ is defined to be $\infty$ if such a maximum does not exist.
\end{notation}

\begin{remark}
	In this article we always assume that the signature of any structure contains the symbol ``$=$'' (represented as equality). This means that the value of $m(\fa)$ is always at least 2.
\end{remark}

\begin{definition}
	Let $\fa$ and $\fb$ be structures. We say that $\fb$ is a \emph{reduct} of $\fa$ iff $\dom(\fa)=\dom(\fb)$ and every relation, function and constant of $\fb$ is first-order definable in $\fa$. Two structures are called \emph{interdefinable} if they are reducts of one another.

	The structures $\fa$ and $\fb$ are called \emph{bidefinable} iff there exist a reduct $\fa'$ of $\fa$ with $\sign(\fa')=\sign(\fb)$, a reduct $\fb'$ of $\fb$ with $\sign(\fb')=\sign(\fa)$, and bijection $i$ from $A$ to $B$ such that $i$ defines an isomorphism from $\fa$ to $\fb'$ and from $\fa'$ to $\fb$. 
\end{definition}

\begin{lemma}\label{interdef}
	The structures $\fa$ and $\fb$ are bidefinable if and only $\fb$ is interdefinable with a structure which is isomorphic to $\fa$.
\end{lemma}

\begin{proof}
	Let us assume that $\fa$ and $\fb$ are bidefinable, and let $\fb'$, $\fa'$ be structures and $i\colon A\rightarrow B$ be a bijection witnessing this. Then by definition $\fa$ and $\fb'$ are isomorphic, and we claim that $\fb$ and $\fb'$ are interdefinable. To see this, it suffices to show that $\fb$ is (isomorphic to) a reduct of $\fb'$. But this follows from the fact that $\fa'$ is a reduct of $\fa$, and $i$ defines an isomorphism from $\fa'$ to $\fb$ and from $\fa$ to $\fb'$. This finishes the proof of the ``only if'' direction of the lemma. 

	Now let us assume that $\fb$ and $\fb'$ are interdefinable, and $\fb'$ and $\fa$ are isomorphic. Let $i$ be an isomorphism from $\fa$ to $\fb'$. Let us define the structure $\fa'$ as follows. The domain set of $\fa'$ is $A$, $\sign(\fa')=\sign(\fb)$, and every symbol in $\sign(\fb)$ is relaized in $\fa'$ is such a way that $i$ defines an isomorphism from $\fa'$ to $\fb$. By definition $\fb'$ is a reduct of $\fb$. It remains to show that $\fa'$ is also a reduct $\fa$. By definition $\fb$ is also a reduct of $\fb'$. Since $i^{-1}$ defines an isomorphism from $\fb'$ to $\fa$ and from $\fb$ to $\fa'$, we obtain that $i^{-1}(\fb)=\fa'$ is a reduct of $i^{-1}(\fb')=\fa$. This finishes the proof that $\fa$ and $\fb$ are bidefinable.
\end{proof}

	From now on we use the description given in Lemma~\ref{interdef} for bidefinability instead of its definition. The reason why we did not define bidefinability this way is that from this description it is not obvious that bidefinability is a symmetric relation.

\begin{remark}
	Clearly interdefinability implies bidefinability, but the converse implication does not hold necessarily even if the structures are assumed to have the same domain set. For instance the structures $(\mathbb{N};=,0)$ and $(\mathbb{N};=,1)$ are bidefinable but not interdefinable.
\end{remark}

	When we are talking about reducts of a given structure $\fa$ we usually consider two reducts the same if they are interdefinable. It turns out that under the assumption that $\fa$ is $\omega$-categorical the reducts of $\fa$ are completely described by the closed supergroups of the automorphism of $\fa$ (see Theorem~\ref{reducts_groups}).

\begin{definition}
	A countable structure $\fa$ is called \emph{$\omega$-categorical} if and only if it is isomorphic to every countable model of its first-order theory.
\end{definition}

	By the theorem of Engeler, Ryll-Nardzewski, and Svenonius, we know that a countably structure $\fa$ is $\omega$-categorical if and only if $\aut(\fa)$ is oligomorphic (see for instance~\cite{HodgesLong}). Given a structure $\fa$ we use the notation $f_n(\fa)$ for $f_n(\aut(\fa))$ where $f$ is any of the operators $\orb,\oi$ or $\os$ (see Definition~\ref{def:growth}). Then Observation~\ref{fewer_orbits} has the following analog for structures and reducts.

\begin{lemma}\label{fewer_orbits_struct}
	Let $\fb$ be a reduct of a structure $\fa$. Then $f_n(\fb)\leq f_n(\fa)$ where $f$ is any of the operators $\orb,\oi$ or $\os$. In particular every reduct of an $\omega$-categorical structure is $\omega$-categorical.
\end{lemma}

	We also know that for an $\omega$-categorical structure $\fa$ a relation (or a function or a constant) defined on $A$ is first-order definable if and only if it is preserved by every automorphism of $\fa$ (see~\cite{Hodges}). This implies the following theorem.

\begin{theorem}\label{reducts_groups}
	Let $\fa$ be a countable $\omega$-categorical structure. Then every closed group $\aut(\fa)\subset G\subset \sym(A)$ is an automorphism group of a reduct of $\fa$, and two reducts of $\fa$ have the same automorphism group if and only if they are interdefinable.
\end{theorem}

	Theorem~\ref{reducts_groups} tells us that given a countable $\omega$-categorical structure $\fa$ the map $\fb\mapsto \aut(\fb)$ defines a bijection between the reducts of $\fa$ (up to interdefinability), and the closed supergroups of $\aut(\fa)$, that is the closed permutation groups containing $\aut(\fa)$.

	An easy consequence of Theorem~\ref{reducts_groups} is the following.

\begin{lemma}\label{bidef_groups}
	Let $\fa$ and $\fb$ be countable structures, so that at least one of them is $\omega$-categorical. Then $\fa$ and $\fb$ are bidefinable if and only $\aut(\fa)$ and $\aut(\fb)$ are isomorphic as permutation groups.
\end{lemma}

\begin{proof}
	We can assume without loss of generality that $\fa$ is $\omega$-categorical.

	Suppose first that $\fa$ and $\fb$ are bidefinable. Then by Lemma~\ref{interdef} we know that $\fb$ is interdefinable with some structure $\fb'$ which is isomorphic to $\fa$. Then $\aut(\fb)=\aut(\fb')$ and it is isomorphic to $\aut(\fa)$.

	For the other direction let $i$ be a bijection from $A$ to $B$ such that $i$ induces an isomorphism from $\aut(\fa)$ to $\aut(\fb)$. Let $\fa':=i^{-1}(\fb)$. Then we have $\aut(\fa')=\iota(i^{-1})\aut(\fb)=\iota(i^{-1})(\iota(i)(\aut(\fa)))=\aut(\fa)$. Theorem~\ref{reducts_groups} then implies that $\fa$ and $\fa'$ are interdefinable. On the other hand it follows from the definition that $\fa'$ is isomorphic to $\fb$. Therefore by Lemma~\ref{interdef} the structures $\fa$ and $\fb$ are bidefinable.
\end{proof}

	Lemma~\ref{bidef_groups}, combined with Theorem~\ref{reducts_groups} implies that for a countable $\omega$-categorical structure the map $\fb\mapsto \aut(\fb)$ defines a bijection between the reducts of $\fa$ up to bidefinability, and the closed supergroups of $\aut(\fa)$ up to isomorphism.

\subsection{Homogeneous structures, Thomas' conjecture}\label{sect:hom}

\begin{definition}
	A structure $\fa$ is called \emph{homogeneous} if every isomorphism between finitely generated substructures of $\fa$ can be extended to an automorphism of $\fa$.
\end{definition}

	We will say that $\fa$ is \emph{finitely related} to abbreviate that $\sign(\fa)$ is finite, and includes no function symbols.

\begin{fact}\label{finrel_ocat}
	Every finitely related homogeneous structure is $\omega$-categorical.
\end{fact}

\begin{proof}
	Let $\fa$ be a finitely related homogeneous structure with signature $\tau$. We show that $\aut(\fa)$ is oligomorphic. Since $\tau$ is finite, and only contains relational and constant symbols, it follows that every finitely generated substructure is finite. Then using the finiteness of $\tau$ again it follows that for all $n\in \omega$ there exist only finitely many $\tau$-structures up to isomorphism. In particular $\fa$ has finitely many quantifier-free types. It follows from the homogeneity of $\fa$ that if two $n$-tuples have the same quantifier-free type then they are in the same orbit of $\aut(\fa)$. This implies that $\orb_n(\fa)$ is finite. Therefore $\aut(\fa)$ is oligomorphic.
\end{proof}

	Simon Thomas conjectured that every finitely related homogeneous structure $\fa$ has only finitely many reducts up to interdefinability~\cite{RandomReducts}.

	The fact that a given $\omega$-categorical structure has finitely many reducts can be formulated in many different ways, as the following proposition shows.

\begin{proposition}\label{thomas_equivalent}
	Let $\fa$ be a countable $\omega$-categorical structure. Then the following are equivalent.
\begin{enumerate}
\item $\fa$ has finitely many reducts up to interdefinability.
\item $\fa$ has finitely many reducts up to bidefinability.
\item There exists some $N\in \omega$ such that every reduct of $\fa$ is interdefinable with a relational structure $\fb$ with $m(\fb)\leq N$.
\item $\aut(\fa)$ has finitely many closed supergroups.
\item $\aut(\fa)$ has finitely many closed non-isomorphic supergroups.
\end{enumerate}
\end{proposition}

\begin{proof}
	The equivalence of items (1) and (2) is Proposition 6.36 in~\cite{bodirsky2018structures}. For the direction (1)$\rightarrow$(3) we only need to show that item (1) implies that every reduct of $\fa$ is interdefinable with a structure which is finitely related. This follows for instance from Lemma 6.35 in~\cite{bodirsky2018structures}.

	In order to show the implication (3)$\rightarrow$(1) it is enough to show that $\fa$ has finitely many reducts $\fb$ with $m(\fb)\leq N$ up to interdefinability. Let $R$ be a relation of $\fb$ with arity at most $N$. Then $R$ is interdefinable with the relation $\{(x_1,\dots,x_N)\colon (x_1,\dots,x_{\ar(R)})\in R\}$. Therefore we can assume without loss of generality that every relation of $\fb$ is $N$-ary. Since $\fb$ is a reduct of $\fa$ it follows that every $N$-ary relation of $\fb$ is a union of some orbits of $\fa$. Therefore there is at most $2^{\orb_N(\fa)}$ many possible choices for an $N$-ary relation of $\fb$, and hence there is at most $2^{2^{\orb_N(\fa)}}$ many choices for reducts $\fb$ of $\fa$ with $m(\fb)\leq N$. Since $\fa$ is $\omega$-categorical we know that $\orb_N(\fa)$ is finite, and hence $2^{2^{\orb_N(\fa)}}$ is also finite.

	The equivalences (1)$\leftrightarrow$(4), and (2)$\leftrightarrow$(5) follow from Theorem~\ref{reducts_groups} and Lemma~\ref{bidef_groups}.
\end{proof}

\begin{definition}
	We call a structure \emph{finitely homogenizable} if it is interdefinable with a homogeneous finitely related structure.

	We denote by $\fhom$ the class of finitely homogenizable structures.
\end{definition}

	Clearly the class $\fhom$ is closed under bidefinability, and by Fact~\ref{finrel_ocat} we know that every structure in $\fhom$ is $\omega$-categorical. Thus it follows from Lemma~\ref{bidef_groups} that the structures in $\fhom$ can fully be described in terms of their automorphisms groups. We give a concrete such description below, as well as other equivalent conditions to check whether a structure is contained in the class $\fhom$. We first need a couple of notations.

\begin{notation}
	For a tuple $\vec{a}\in X^n$ we denote by $a_i$ the $i$th coordinate of $\vec{a}$. For a function $\pi\colon k\rightarrow n$ we denote by $\vec{a}_{\pi}$ the tuple $(a_{\pi(0)},\dots,a_{\pi(k-1)})$.
\end{notation}

\begin{notation}
	Let $\fa$ be a relational structure, and $m\in \omega$. Then we denote by $\typestruct_m(\fa)$ the structure whose domain is $A$, and whose relations are all subsets of $A^m$ which are definable in $\fa$.
\end{notation}

	It is clear from the definition that if $m'\leq m$ then every $m'$-ary relation of $\fa$ is quantifier-free definable in $\typestruct_m(\fa)$, that is $\typestruct_{m'}(\fa)$ is a quantifier-free reduct of $\typestruct_m(\fa)$.

	If $\fa$ is $\omega$-categorical then every type over $\fa$ is principal, and thus in this case $\typestruct_m(\fa)$ is a first-order reduct of $\fa$ for all $m\in \omega$.

	For a structure $\fa$ and tuples $\vec{u},\vec{v}\in A^n$ we write $\vec{u}\same_{\fa}\vec{v}$ iff $\vec{u}$ and $\vec{v}$ have the same type over $\fa$, and we write $\vec{u}\notsame_{\fa} \vec{v}$ iff $\vec{u}\same_{\fa} \vec{v}$ does not hold. In the case when $\fa$ is $\omega$-categorical $\vec{u}\same_{\fa} \vec{v}$ holds if and only if $\vec{u}$ and $\vec{v}$ are in the same orbit of $\aut(\fa)$.

\begin{lemma}\label{when_finhom}
	Let $\fa$ be an $\omega$-categorical structure, and let $m\geq 2$. Then the following are equivalent.
\begin{enumerate}
\item $\fa$ is interdefinable with a homogeneous relational signature structure $\fa'$ with $m(\fa')\leq m$.
\item $\fa$ is interdefinable with a homogeneous relational signature structure all of whose relations have arity exactly $m$.
\item $\fa$ and $\typestruct_m(\fa)$ are interdefinable, and $\typestruct_m(\fa)$ is homogeneous.
\item For all $\vec{u},\vec{v}\in A^n$ if $\vec{u}\notsame_{\fa}\vec{v}$ then there exists a function $\pi\colon m\rightarrow n$ such that $\vec{u}_{\pi}\notsame_{\fa}\vec{v}_{\pi}$.
\end{enumerate}
\end{lemma}

\begin{proof}
	The implication (2)$\rightarrow$(1) is obvious. The converse follows from the fact that if $m'\leq m$ then every $m'$-ary relation of $\fa'$ is quantifier-free definable from a $m$-ary relation of $\fa'$.

	The implication (3)$\rightarrow$(2) is obvious. For the converse implication let us assume that item (2) holds, and let $\fa'$ be a structure as it is dictated in the condition. Then $\typestruct_m(\fa)=\typestruct_m(\fa')$. Since $\fa$, and thus also $\fa'$ is $\omega$-categorical it follows that every definable relation of $\fa'$ of arity $m$ is a finite union of some $m$-orbits of $\fa'$. Therefore $\fa'$ is quantifier-free definable in $\typestruct_m(\fa')$. Since $\typestruct_m(\fa')$ is reduct of $\fa'$ and $\fa'$ is homogeneous it follows that $\typestruct_m(\fa')$ is also quantifier-free definable in $\fa'$. Therefore $\typestruct_m(\fa')$ is interdefinable with $\fa$, and it is homogeneous.

	(3)$\rightarrow$(4). Let us assume that condition (3) holds. We can assume without loss of generality that $\fa=\typestruct_m(\fa)$. We show the contrapositive of the statement in item (4). Let us assume that for all $\pi\colon m\rightarrow n$ we have $\vec{u}_{\pi}\notsame_{\fa}\vec{v}_{\pi}$. Then by definition we have $R(\vec{u}_{\pi})\Leftrightarrow R(\vec{v}_{\pi})$ for all $\pi\colon m\rightarrow n$ and $\ar(R)=m$. By definition the arity of every relation of $\fa$ is exactly $m$. This means that the tuples $\vec{u}$ and $\vec{v}$ have the same atomic type over $\fa$. Since $\fa$ is homogeneous, it has quantifier elimination, and hence $\vec{u}$ and $\vec{v}$ have the same type over $\fa$.

	(4)$\rightarrow$(3). We show the contrapositive. Let us assume that item (3) does not hold. This means that $\fa$ has a relation $R$ which is not quantifier-free definable over $\typestruct_m(\fa)$. Then there exist tuples $\vec{u}\in R, \vec{v}\in \bar{R}$ such that for every relation $S$ of $\typestruct_m(\fa)$, and $\pi\colon m\rightarrow n$ we have $\vec{u}_{\pi}\in S$ if and only if $\vec{v}_{\pi}\in S$. By definition this implies that for all $\pi\colon m\rightarrow n$ the tuples $\vec{u}_{\pi}$ and $\vec{v}_{\pi}$ have the same type over $\fa$. Thus the tuples $\vec{u}$ and $\vec{v}$ witness that Condition (4) fails.
\end{proof}

	The following is a direct consequence of Lemma~\ref{when_finhom}.

\begin{corollary}\label{when_finhom2}
	Let $\fa$ be an $\omega$-categorical structure. Then the following are equivalent.
\begin{enumerate}
\item $\fa\in \fhom$.
\item $\fa$ is interdefinable with $\typestruct_m(\fa)$ for some $m$.
\item There exists an $m\in \omega$ such that for all $\vec{u},\vec{v}\in A^n$ if $\vec{u}\notsame_{\fa}\vec{v}$ then there exists a function $\pi\colon m\rightarrow n$ such that $\vec{u}_{\pi}\notsame_{\fa}\vec{v}_{\pi}$.
\end{enumerate}
\end{corollary}

\subsection{Finite boundedness}\label{sect:fin_bound}

	In the study of finitely related homogeneous structures we are particularly interested in the ones which are \emph{finitely bounded}.

\begin{definition}
	Let $\tau$ be a relational signature, and let $\mathcal{C}$ be a set of finite $\tau$-structures. Then we denote by $\forb(\mathcal{C})$ the class of those finite $\tau$-structures which do not embed any structure from $\mathcal{C}$.

	A relational structure $\fa$ is called \emph{finitely bounded} if its signature $\tau:=\sign(\fa)$ is finite, and there is finite set $\mathcal{C}$ of $\tau$-structures such that $\age(\fa)=\forb(\mathcal{C})$.
\end{definition}

	Next we give a couple of equivalent conditions for a finitely related homogeneous structure to be finitely bounded which we will need later.

\begin{notation}
	Let $\mathcal{C}$ be a class of finite $\tau$-structures. Then 
\begin{itemize}
\item we denote by $\comp{\mathcal{C}}$ the class of those finite $\tau$-structures which are not in $\mathcal{C}$. 
\item we denote by $\Min(\mathcal{C})$ the class of minimal structures in $\mathcal{C}$ i.e., the class of those structures in $\mathcal{C}$ which do not have a nontrivial substructure which is also in $\mathcal{C}$.
\end{itemize}
\end{notation}

\begin{lemma}\label{when_finbounded}
	Let $\fa$ be a finitely related structure with signature $\tau$. Then the following are equivalent.
\begin{enumerate}
\item $\fa$ is finitely bounded.
\item $\Min(\comp{\age(\fa)})$ contains finitely many structures up to isomorphism.
\item There exists some $m$ such that the size of every structure in $\Min(\comp{\age(\fa)})$ is at most $m$.
\end{enumerate}
\end{lemma}

\begin{proof}
	The equivalence of conditions (2) and (3) is clear since $\tau$ is finite.

	(2)$\rightarrow$(1). Let $\mathcal{F}$ be a finite subset of $\Min(\comp{\age(\fa)})$ so that every structure in $\Min(\comp{\age(\fa)})$ is represented in $\mathcal{F}$. We claim that $\age(\fa)=\forb(\mathcal{F})$. Let us suppose that $\fb\not\in\forb(\mathcal{F})$. Then $\fb$ has a substructure $\fb'$ which is isomorphic to some structure in $\mathcal{F}$. In particular $\fb'\not\in \age(\fa)$. Since $\age(\fa)$ is closed under taking substructures it follows that $\fb\not\in \age(\fa)$. For the other direction let us assume that $\fb\not\in \age(\fa)$. Then $\fb$ contains a minimal structure with this property, say $\fb'$. Then by definition $\fb'$ is isomorphic to some structure in $\fc\in \mathcal{F}$, and hence $\fc$ embeds into $\fb$, that is $\fb\not\in \forb(\mathcal{F})$.

	(1)$\rightarrow$(3). Let us assume that $\age(\fa)=\forb(\mathcal{F})$ for some finite set $\mathcal{F}$ of $\tau$-structures, and let $m$ be the maximum of the sizes of structures in $\mathcal{F}$. Let $\fb\in \Min(\comp{\age(\fa)})$. Since $\fb\not\in \age(\fa)=\forb(\mathcal{F})$ it follows $\fb$ has a substructure $\fb'$ which is isomorphic to some $\fc\in \mathcal{F}$. But in this case $\fb'$ is also not in $\age(\fa)$. Then the minimality of $\fb$ implies that in fact $\fb=\fb'$, and hence $|\fb|\leq m$.
\end{proof}

\begin{notation}
	For a structure $\fa$ we denote by $b(\fa)$ the maximum of $m(\fa)$ and the maximum of sizes of structures in $\Min(\comp{\age(\fa)})$.
\end{notation}

	Lemma~\ref{when_finbounded} shows that if $\fa$ is homogeneous then it is finitely bounded if and only if $b(\fa)$ is finite.

\begin{notation}
	We denote by $\fbh$ the class of those structures which are interdefinable with a finitely bounded homogeneous relational structure.
\end{notation}

	We can argue the same way as in the previous section that the structures $\fbh$ can be described in terms of their automorphisms. In this case however we do not know of any nice such description. Nevertheless we know that the class $\{\aut(\fa)\colon\fa\in \fbh\}$ is closed under taking finite direct products, infinite copies, and rather surprisingly (as we will see in Section~\ref{sect:fin_bounded}) it is also closed under taking finite index supergroups.	

\subsection{Unions and copies of structures}

	Let us consider the following two constructions on structures.

\begin{definition}\label{def:union}
	Let $\fa_1,\dots,\fa_n$ be relational structures with pairwise disjoint domains. 

	Then we define the \emph{disjoint union} $\fb$ of the structures $\fa_1,\dots,\fa_n$, denoted by $\biguplus_{i=1}^n\fa_i$ such that
\begin{itemize}
\item the domain set of $\fb$ is $\bigcup_{i=1}^nA_i$,
\item its signature is $\bigcup_{i=1}^n\sign(\fa_i)$ together with some unary symbols $U_1,\dots,U_n\not\in \bigcup_{i=1}^n\sign(\fa_i)$,
\item $U_i^\fb=A_i$ for all $1\leq i\leq n$ and,
\item $R^\fb=R^{\fa_i}$ for all $R\in \sign(\fa_i)$.
\end{itemize}
\end{definition}

\begin{definition}\label{def:copies}
	Let $\fa$ be a relational structure. Then we define the structure $\fa\wr \omega$, called \emph{infinitely many copies of $\fa$} as follows. 
\begin{itemize}
\item The domain set of $\fa\wr \omega$ is $A\times \omega$, 
\item its signature is $\sign(\fa)\sqcup \{E\}$ where $\ar(E)=2$,
\item $E^{\fa\wr \omega}=\{(a_1,n),(a_2,n)\colon a_1,a_2\in A, n\in \omega\}$,
\item $R^{\fa\wr \omega}=\{(a_1,n),\dots,(a_k,n)\colon (a_1,\dots,a_k)\in R^{\fa}, n\in \omega\}$ for all $R\in \sign(\fa)$,
\end{itemize}
\end{definition}

	The following observation is an easy consequence of the definitions.

\begin{proposition}\label{aut_union}
	Let $\fa_1,\dots,\fa_n$ be relational structures as in Definition~\ref{def:union}. Then $\aut(\biguplus_{i=1}^n\fa_i)=\prod_{i=1}^n\aut(\fa_i)$.
\end{proposition}

\begin{proposition}\label{aut_copies}
	Let $\fa$ be a relational structure. Then $\aut(\fa\wr \omega)=\aut(\fa)\wr \sym(\omega)$.
\end{proposition}

	The following two lemmas plays a crucial role in our proofs later.

\begin{lemma}\label{union_homogeneous}
	Let $\fa_1,\dots,\fa_n$ be structures with pairwise disjoint domains, and let $\fa:=\prod_{i=1}^n\fa_i$. Then
\begin{enumerate}
\item if $\fa_1,\dots,\fa_n$ are all homogeneous, then so is $\fa$, and
\item if $\fa_1,\dots,\fa_n$ are all finitely bounded, then so is $\fa$.
\end{enumerate}
\end{lemma}

\begin{proof}
	Item (1) follows from a straightforward induction using Lemma 3.1 in~\cite{BodirskyRamsey}.

	As for item (2) let us assume that $\age(\fa)=\forb(\mathcal{F}_i)$ for some finite sets of structures $\mathcal{F}_i\colon i=1,\dots,n$. For a $\sign(\fa_i)$-structure $\fb$ let $(\fb;U_i)$ denote the structure $\fb$ expanded by a unary relation denoted by $U_i$ which contains all elements of $B_i$. Let $\mathcal{G}$ be the set of one-element $\{U_1,\dots,U_n\}$-structures $\fb$ such that either $\fb^{U_i}=\emptyset$ for all $i$ or $\fb^{U_i}=\fb^{U_j}=B$ for some $i\neq j$. Then it follows from the construction that $$\age(\fa)=\forb(\{(\fb,U_i)\colon 1\leq i\leq n, \fb\in \mathcal{F}_i\}\cup \mathcal{G}).$$ Therefore $\fa$ is finitely bounded.
\end{proof}

\begin{lemma}\label{copies_homogeneous}
	Let $\fa$ be homogeneous structures. Then $\fa\wr \omega$ is homogeneous, and if $\fa$ is finitely bounded, then so is $\fa\wr \omega$.
\end{lemma}

\begin{proof}
	Let us assume that the tuples $\vec{u}=((a_1,m_1),\dots,(a_k,m_k)),\vec{v}=((b_1,n_1),\dots,(b_k,n_k))$ have the same quantifier-free type over $\fa\wr \omega$. We have to show that $\vec{u}$ and $\vec{v}$ are in the same orbit of $\aut(\fa\wr \omega)=\aut(\fa)\wr \sym(\omega)$. For an $n\in \omega$ let us denote by $\vec{u}_n$ ($\vec{v}_n$) the subtuple of $\vec{u}$ ($\vec{v}$) containing those elemtnens whose second coordinate is $n$. By definition it follows that $m_i=m_j$ iff $n_i=n_j$ for all $1\leq i,j\leq k$. Thus there exists a permutation $\sigma\in \sym(\omega)$ such that $\vec{u}_n:=((a_{i_1},n),\dots,a_{i_k},n))$ we have $n_{i_1}=\dots=n_{i_k}=\sigma(b)$. Moreover by definition it follows that the tuples $(a_{i_1},\dots,a_{i_k})$ and $(b_{i_1},\dots,b_{i_k})$ have the same quantifier-free type, and thus by the homogeneity of $\fa$ it follows that there exists some $\alpha_n\in \aut(\fa)$ such that $\alpha(a_{i_j})=b_{i_j}$ for all $j=1,\dots,k$. Now let $\gamma\in \sym(A\times \omega)$ be defined as $\gamma(a,n):=(\alpha_n(a),\sigma(n))$. Then $\gamma\in \aut(\fa)\wr \sym(\omega)$ and by our construction $\gamma(a_i,m_i)=(\alpha_{m_i}(a_i),\beta(m_i))=(c_i,d_i)$ for all $1\leq i\leq k$.

	For the second part of the lemma let us assume that $\age(\fa)=\forb(\mathcal{F})$ for some finite set of structures $\mathcal{F}$. Let $\mathcal{G}_0$ be a finite set of forbidden structures in the signature $\{E\}$ so that $\fa\in \forb(\mathcal{G}_0)$ if and only if $E^{\fa}$ is an equivalence relation on $\fa$. We can do this by listing all 3-element $\{E\}$-structures where $E$ does not define an equivalence relation. For a $\sign(\fa)$-structure $\fc$ let $\fc^*$ denote the structure $\fc$ expanded by a binary relation, denoted by $E$ which contains every pair in $C^2$. Let $$\mathcal{G}:=\mathcal{G}_0\cup \{\fc^*\colon \fc\in \mathcal{F}\}.$$ Then it is straightforward to check that $\forb(\mathcal{G})$ is exactly the age of $\fa\wr \omega$. Therefore $\fa\wr \omega$ is finitely bounded.
\end{proof}

\begin{corollary}\label{fbh_union_closed}
	The classes $\fhom$ and $\fbh$ are closed under taking finite disjoint unions and infinite copies.
\end{corollary}

\section{Hereditarily cellular structures}\label{sect:class}

\begin{definition}
	Let $T$ be a complete theory is some language $L$. Then we say that $T$ is \emph{monadically stable} if every expansion of $T$ by unary predicates is stable. A first-order structure $\fa$ is called \emph{monadically stable} if its first-order theory is monadically stable.

	We denote the class of countable $\omega$-categorical monadically stable structures by $\ms$. For a natural number $k$ we denote by $\ms_k$ the class of structures in $\ms$ whose Morley rank is at most $k$.
\end{definition}

	Lachlan showed in~\cite{lachlan1992} that the structures in $\ms$ can be described in terms of their automorphism group. For the description of these groups we first need a couple of definitions.

\begin{definition}\label{omega_part}
	Let $G\leq \sym(X)$ be a permutation group. Then a triple $\partition=(K, \bigeq, \smalleq)$ is called an \emph{$\omega$-partition} of $G$ iff
\begin{enumerate}
\item $K\subset X$ is finite, and it is fixed setwise by $G$.
\item $\smalleq$ and $\bigeq$ are congruences of $G|_{X\setminus K}$ with $\smalleq\subset \bigeq$.
\item $\bigeq$ has finitely many classes.
\item Every $\bigeq$-class is a union of $\aleph_0$ many $\smalleq$-classes.
\item For every $C\in X/\bigeq$ we have $G_{((C))}/\smalleq=\sym(C/\smalleq)$.
\end{enumerate}
	The \emph{components} of an $\omega$-partition $\partition=(K, \bigeq, \smalleq)$ are the permutation groups $G_{((Y))}$ for all $Y\in X/\smalleq$.
\end{definition}

\begin{remark}
	Definition~\ref{omega_part} includes the case when $X=K$. In this case $\smalleq=\bigeq=\emptyset$, and therefore all conditions (2)-(5) in Definition~\ref{omega_part} are automatically satisfied.
\end{remark}

\begin{definition}
	We define the classes of closed permutation groups $\gone_{-1},\gone_0,\gone_1,\dots$ recursively as follows.
\begin{itemize}
\item $\gone_{-1}=\emptyset$,
\item For $n\geq 0$ we define $G\in \gone_n$ if and only if $G$ has an $\omega$-partition $\partition$ so that every component of $\partition$ is in $\gone_{n-1}$.
\end{itemize}

	We define the class $\gone$ as the union of all classes $\gone_n\colon n\in \omega$. We call $\gone$ the class of \emph{hereditarily $\omega$-partitioned} permutation groups. The \emph{rank} of a hereditarily $\omega$-partitioned group $G$, denoted by $\rk(G)$ is the smallest $k$ such that $G\in \gone_n$.
	
	We say that a structure $\fa$ is \emph{hereditarily cellular} if and only if $\aut(\fa)$ hereditarily $\omega$-partitioned.
\end{definition}

\begin{remark}\label{rank0}
	It is clear from the definition that $\gone_0$ contains exactly those groups which act on a finite set. Furthermore it can easily by shown by induction that $\gone_{n-1}\subset \gone_n$ for all $n\in \omega$.
\end{remark}

\begin{theorem}[\cite{lachlan1992}, Theorem 2.3]\label{lachlan}
	$\fa\in \ms$ if and only if $\fa$ is hereditarily cellular.
\end{theorem}

	There are two more facts from~\cite{lachlan1992} that we need to keep in mind. We quote them below, slightly rephrased.

\begin{proposition}[\cite{lachlan1992}, Proposition 1.4]\label{prop:lachlan14}
	Let $(K, \bigeq, \smalleq)$ and $(K', \bigeq', \smalleq')$ be $\omega$-partitions of a permutation group $G\in\gone$. Then $K=K'$, $\bigeq=\bigeq'$ and $\smalleq=\smalleq'$.
\end{proposition}

\begin{proposition}[\cite{lachlan1992}, Propositions 1.3 and 1.6]\label{prop:lachlan16}
	Let $\fa\in \ms$. Then $\fa$ is $\omega$-stable, and its Morley rank is equal to $\rk(\aut(\fa))$.
\end{proposition}

	By Proposition~\ref{prop:lachlan14}, the $\omega$-partition witnessing $G\in\gone$ is unique. Therefore from now on we will refer to this $\omega$-partition as \emph{the} $\omega$-partition of the group $G$. Although our definition of the rank function $\rk$ on $\gone$ is different from the one given in~\cite{lachlan1992} it is easy to see that the two definitions are equivalent. Then by Proposition 1.4 in~\cite{lachlan1992} we know that if $\fa\in \ms$ then $\rk(\aut(\fa))$ is equal to the Morley rank of $\fa$.

\begin{corollary}\label{lachlan2}
	Let $\fa$ be a structure. Then $\fa\in \ms_n$ if and only if $\aut(\fa)\in \gone_n$. In particular $\ms=\bigcup_{n=0}^\infty \ms_n$, that is every structure in $\ms$ has a finite Morley rank.
\end{corollary}

\begin{example}\label{def:single}
	We define the structures $\single_n$ by recursion as follows. $\single_0$ is defined to be the 1-element pure set, and $\single_n=\single_{n-1}\wr \omega$. The structure $\single_n$ can also be described as a countable infinite set $A$ endowed with equivalence classes $E_0\subsetneq \dots\subsetneq E_n$ such that
\begin{itemize}
\item $E_0$ is the trivial equivalence relation,
\item $E_n$ is the universal equivalence relation, and
\item every $E_i$ class contains $\aleph_0$ many $E_{i-1}$ classes for all $1\leq i\leq n$.
\end{itemize}

	It is easy to see that $\single_n\in \ms$, and $\rk(\single_n)=n$.	
\end{example}

	The class $\ms$ can also be described by orbit growth on subsets. We say that a function $f\colon\omega\rightarrow \mathbb{R}$ is \emph{slower (faster) then exponential} if for all $c>1$ we have $f(n)<c^n$ ($f(n)>c^n$) if $n$ is large enough. Using these notions we know the following dichotomy for $\omega$-categorical stable structures.

\begin{theorem}[\cite{braunfeld2019monadic}, Theorem 3.10]\label{sam_original}
	Let $\fa$ be an $\omega$-categorical stable structure. Then exactly one of the following holds.
\begin{enumerate}
\item $\os_n(\fa)$ is slower than exponential, and $\fa$ is monadically stable.
\item $\os_n(\fa)$ is faster than exponential, and $\fa$ is not monadically stable.
\end{enumerate}
\end{theorem}

	The following is a straightforward consequence of Theorem~\ref{sam_original}.
	
\begin{theorem}\label{sam}
	Let $\fa$ be a structure. Then $\fa\in \ms$ if and only if $\fa$ is stable, and $\os_n(\fa)$ is slower than exponential.
\end{theorem}

	Note that if $\fb$ is a reduct of $\fa$ then $\os_n(\fb)\leq \os_n(\fa)$ for all $n$. Therefore Theorem~\ref{sam} immediately implies that the class $\ms$ is closed under taking reducts. In fact one can arrive to the same conclusion directly from the definition. Moreover if $\fb$ is a reduct of $\fa$ then every relation which is definable over an elementary expansion of $\fb$ is also definable over an elementary expansion of $\fa$. Hence the Morley rank of $\fb$ is at most the Morley rank of $\fa$. This implies in turn that for all $n\in\omega$ the class $\ms_n$ is also closed under taking reducts. In terms of permutation groups, this means that the classes $\gone_n\colon n\in \omega$, and $\gone$ are closed under taking closed supergroups (see Theorem~\ref{reducts_groups} and Corollary~\ref{lachlan2}).
 
	We can also use Theorem~\ref{sam} to show that the class $\ms$ is closed under taking model-complete cores. As explained in the introduction this observation is particularly important in the study of CSPs.
	
\begin{definition}
	An $\omega$-categorical structure is called a \emph{model-complete core} if the topological closure of $\aut(\fa)$ is equal to the endomorphism monoid of $\fa$.
	
	We say that two structures $\fa$ and $\fb$ are \emph{homomorphically equivalent} if there exist homomorphisms $f\colon\fa\rightarrow \fb$ and $g\colon\fb\rightarrow \fa$.
\end{definition}

\begin{theorem}[~\cite{Cores-journal}, Theorem 16]\label{mc_core}
	Let $\fa$ be an $\omega$-categorical structure. Then $\fa$ is homomorphically equivalent to a model-complete core $\fb$. Moreover the structure $\fb$ is unique up to isomorphism, and it is again $\omega$-categorical.
\end{theorem}

	We call the (up to isomorphism) unique structure $\fb$ as in Theorem~\ref{mc_core} \emph{the} model-complete core of $\fa$.
	
\begin{lemma}\label{mc_orbit}
	Let $\fa$ be an $\omega$-categorical structure, and let $\fb$ be its model-complete core. Then $f_n(\fb)\leq f_n(\fa)$ for all $n\in \omega$ where $f$ is any of the operators $\orb,\oi$ or $\os$.
\end{lemma}

\begin{proof}
	The case when $f=\orb$ is shown in~\cite{Bodirsky-HDR}. The statement for $f=\oi$ and $f=\os$ can be shown analogously.
\end{proof}

	Combining Theorem~\ref{sam} and Lemma~\ref{mc_orbit} we obtain that the class $\ms$ is closed under taking model-complete cores. By refining this argument one can also show that the same conclusion holds for all classes $\ms_n\colon n\in \omega$. 

\subsection{Building up the class $\gone$.}\label{sect:build}

	In this section we construct recursively a subclass $\gone'$ of $\gone$, and examine some properties of groups in $\gone'$. Then in the next section we will show that in fact $\gone'=\gone$, and as a result we obtain an explicit recursive description of $\gone$.

\begin{lemma}\label{inf_copies}
	Let $H$ be a permutation group. Then the group $G:=H\wr \sym(\omega)$ has an $\omega$-partition with all components isomorphic to $H$.
\end{lemma}

\begin{proof}
	Let $\smalleq:=(\dom(H)\times \omega)^2$ and $\bigeq:=\{((a,n),(b,n))| a,b\in H, n\in \omega\}$. Then $\partition:=(\emptyset,\smalleq,\bigeq)$ defines an $\omega$-partition of
$G$ with all components isomorphic to $H$.
\end{proof}

\begin{corollary}\label{inf_copies2}
	If $H\in \gone_n$ then $H\wr \sym(\omega)\in \gone_{n+1}$.
\end{corollary}

	The following lemma tells us that the classes $\gone_n$ are closed under finite direct products.

\begin{lemma}\label{direct_gone}
	Let $H_1,\dots,H_k\in \gone_n$ with pairwise disjoint domains. Then $H:=\prod_{i=1}^kH_i\in \gone_n$.
\end{lemma}

\begin{proof}
	Let $\partition_i:=(K_i,\bigeq_i,\smalleq_i)$ be the $\omega$-partition of $H_i$ witnessing the fact that $H_i\in \gone_n$, that is every component of $\partition_i$ is in $\gone_{n-1}$. Then let $\partition:=(K, \bigeq, \smalleq)$ with $K:=\bigcup_{i=1}^kK_i, \bigeq:=\bigcup_{i=1}^k(\bigeq_i)$ and $\smalleq:=\bigcup_{i=1}^k(\smalleq_i)$. We claim that $\partition$ is an $\omega$-partition of $H$, and every component of $\partition$ is a component of one of the $\omega$-partitions $\partition_i: i=1,\dots,k$. Since these components are all contained in $\gone_{n-1}$ by our assumption, this implies the statement of the lemma.

	It follows from the definition that $H$ fixes the sets $X_i:=\dom(H_i)$, and $H|_{X_i}=H_i$. This implies immediately that $G$ fixes the set $K$, and the equivalence relations $\smalleq$ and $\bigeq$. It is clear from the definition that $K$ is finite, $\bigeq$ has finitely many classes, and every $\bigeq$ class is a union of $\aleph_0$ many $\smalleq$ classes. It remains to show that item (5) of Definition~\ref{omega_part} holds for the congruences $\bigeq$ and $\smalleq$. Let $C$ be a $\bigeq$ class. Then by definition $C$ is also a $\bigeq_i$ class for some $1\leq i\leq k$. Then we have 
\begin{multline*}
H_{((C))}/\smalleq=H_{((C))}/\smalleq_i=(H|_{X_i})_{((C))}/\smalleq_i=\\=(H_i)_{((C))}/\smalleq_i=\sym(C/\smalleq_i)=\sym(C/\smalleq),
\end{multline*}
	where the second to last equality above follows from the assumption that $\partition_i$ is an $\omega$-partition. This finishes the proof that $\partition=(K,\bigeq,\smalleq)$ is an $\omega$-partition. Now if $X$ is a $\smalleq$ class then $X$ is also a $\smalleq_i$ class for some $1\leq i\leq k$, and thus $H_{((Y))}=(H|_{X_i})_{((Y))}=(H_i)_{((Y))}$ which is by definition a component of $\partition_i$. This shows that every components of $\partition$ is a component of $\partition_i$ for some $1\leq i \leq k$.
\end{proof}

	In the next few paragraphs, we will construct a structure with an $\omega$-partition and dissect a permutation that preserves the equivalence relations on it. Our ultimate goal is to give an explicit description of the groups in $\gone$ and the structures in $\ms$.

	Now we fix some notational conventions that we will use throughout this section. First we fix some pairwise disjoint sets $Y_0,\ldots,Y_k$, with $Y_0$ finite and all other $Y_i$ nonempty; let $Y:=\bigcup_{i=0}^kY_i$. We define $i\colon Y\to\{0,\ldots,k\}$ as the function mapping $a$ to the unique $j$ such that $a\in Y_j$. We will use $\Omega(Y_0,\dots,Y_k)$ to denote the set $\bigcup_{i=1}^k(Y_i\times \omega)\cup Y_0\times \{0\}$. Given $u=(a,n)\in \Omega(Y_0,\dots,Y_k)$, we define $p(u)=a, q(u)=n$.

Let $$\bigeq:=\bigcup_{i=1}^k(Y_i\times \omega)^2\cup (Y_0\times\{0\})^2$$ and $$\smalleq:=\bigcup_{i=1}^k(\{(a,n),(b,n)| a,b\in Y_i, n\in \omega\})\cup (Y_0\times\{0\})^2.$$  Observe that $Y_0\times\{0\}$ is the only finite $\bigeq$-class, so any permutation of $\Omega(Y_0,\dots,Y_k)$ that preserves $\bigeq$ must stabilize $Y_0\times\{0\}$ setwise.

Suppose that $\sigma\in\sym(\Omega(Y_0,\dots,Y_k))$ preserves $\bigeq$ and $\smalleq$. We identify $\Omega(Y_0,\dots,Y_k)/\bigeq$ with $\{0,\ldots,k\}$ via the mapping that sends $Y_0\times\{0\}$ to 0 and $Y_i\times\omega$ to $i$, for $1\leq i\leq k$. Similarly, $\Omega(Y_0,\dots,Y_k)/\smalleq$ is identified with $\{(0,0)\}\cup\{1,\ldots,k\}\times\omega$. These identifications allow us to think of $\sigma/\bigeq$ as a permutation $\phi(\sigma)\in\sym(\{0,\ldots,k\})$ that fixes 0 and similarly $\sigma/\smalleq$ as a permutation $\psi(\sigma)\in\sym(\{(0,0)\}\cup\{1,\ldots,k\}\times\omega)$ that fixes $(0,0)$. 

We can further decompose $\psi(\sigma)$ into a constant mapping $\psi_0\equiv 0$ and $k$ permutations of $\omega$, say $\psi_i(\sigma)$, so that $\psi(\sigma)(i,n)=(\phi(\sigma)(i),\psi_i(\sigma)(n))$ for all $(n,i)\in \{(0,0)\}\cup (\{1,\dots,k\}\times\omega)$. Thus $\psi(\sigma)$ is uniquely determined by the $\psi_i(\sigma)$ ($0\leq i \leq k$) and $\phi(\sigma)$. 

Given any $\vec{n}=(0,n_1,\ldots,n_k)\in\{0\}\times\omega^k$, we define $\vec{\mathbf{n}}$ as $((0,0),(1,n_1),\ldots,(k,n_k))\in\{(0,0)\}\cup\{1,\ldots,k\}\times\omega$. For each $\vec{n}$ we define an injection $\ell_{\vec{n}}\colon Y\to\Omega(Y_0,\ldots,Y_k)$, defined by $\ell_{\vec{n}}(a)=(a,n_{i(a)})$. We define $$\rho_{\vec{n}}(\sigma)\coloneqq p\circ\sigma\circ\ell_{\vec{n}}.$$ These functions are elements of $\sym(Y)$. For $n\in\omega$ and $i\in\{0,\ldots k\}$ we define $$\rho_n^i(\sigma)\coloneqq\rho_{\mathbf{n}}(\sigma)|_{Y_i},$$ where $\mathbf{n}$ is the tuple $(0,n,\ldots,n)$. Then $\rho_n^i(\sigma)$ is a bijection $Y_i\to Y_{\phi(\sigma)(i)}$.

	Let $\sigma\in \aut(\Omega(Y_0,\dots,Y_k);Y_0\times\{0\},\bigeq,\smalleq)$ be arbitrary. Then using the notation above it follows that given $(a,n)\in \Omega(Y_0,\dots,Y_k)$ we have $\sigma((a,n))=(\rho_n^{i(a)}(\sigma)(a),\psi_{i(a)}(\sigma)(n))$. In particular a permutation $\sigma$ in $\aut(\Omega(Y_0,\dots,Y_k);Y_0,\bigeq,\smalleq)$ is uniquely determined by the maps $\rho_{n}^i(\sigma)$ and $\psi_i(\sigma)$. 

\begin{lemma}\label{rho_product}
	Let $\sigma,\tau\in\aut(\Omega(Y_0,\dots,Y_k);Y_0\times\{0\},\bigeq,\smalleq)$. Then for all $\vec{n}=(n_0,\dots,n_k)\in \{0\}\times \omega^k$ we have
\begin{enumerate}[(i)]
\item $\rho_{\vec{n}}(\tau\sigma)=\rho_{\psi(\sigma)(\vec{\mathbf{n}})}(\tau)\rho_{\vec{n}}(\sigma)$, and \item $\rho_{\vec{n}}(\sigma^{-1})=(\rho_{\psi(\sigma^{-1})(\vec{\mathbf{n}})}(\sigma))^{-1}$.
\end{enumerate}
\end{lemma}

\begin{proof}
	(i) Let $a\in Y_i$. Then $\sigma$ maps $u=(a,n_i)$ to $u':=(\rho_{\vec{n}}(\sigma)(a),\psi_{i'}(\vec{n}))$ for $i'=i(\rho_{\vec{n}}(\sigma)(a))$. By definition $p(\tau(u'))=\rho_{\psi(\sigma)(\vec{\mathbf{n}})}(\tau)(\rho_{\vec{n}}(\sigma)(a))$. On the other hand $p(\tau(u'))=p(\tau\sigma(u))=\rho_{\vec{n}}(\tau\sigma)(a)$. This shows the equality in item (i).

	(ii) Using item (i) it follows that $$\id(Y)=\rho_{\vec{n}}(\id)=\rho_{\vec{n}}(\sigma\sigma^{-1})=\rho_{\psi(\sigma^{-1})(\vec{\mathbf{n}})}(\sigma)\rho_{\vec{n}}(\sigma^{-1}),$$ and the equality (ii) follows.
\end{proof}

Suppose now that $N_0,\dots,N_k$ are permutation groups with $\dom(N_i)=Y_i$, and $N_0=\id(Y_0)$, the trivial permutation group on $Y_0$. We will use the notation $N$ for the group $\prod_{i=0}^kN_i$.

\begin{definition}\label{def:groups2}
	Let $H$ be a permutation group acting on $Y$ such that $N$ is a normal subgroup of $H$, and $H$ fixes $Y_0$ setwise. Then we define the subset $\group(H;N_0,\dots,N_k)$ of $\sym(\Omega(Y_0,\dots,Y_k))$ as follows. A permutation $\sigma$ is contained in $\group(H;N_0,\dots,N_k)$ if and only if
\begin{enumerate}[(a)]
\item\label{item:y0} $\sigma$ fixes $Y_0\times\{0\}$ setwise,
\item\label{item:eq} $\sigma$ preserves the relations $\bigeq$ and $\smalleq$,
\item\label{item:h} for all $\vec{n}=(n_0,n_1,\dots,n_k)\in \{0\}\times \omega^k$ it holds that $\rho_{\vec{n}}(\sigma)\in H$, and
\item\label{item:mod} for all $\vec{n},\vec{m}\in \{0\}\times \omega^k$ the maps $\rho_{\vec{n}}(\sigma)$ and $\rho_{\vec{m}}(\sigma)$ are equal modulo $N$.
\end{enumerate}

	Define $$\smallgroup(N_0,\dots,N_k):=\id(Y_0\times \{0\})\times \prod_{i=1}^k(N_i\wr \sym(\omega))\subset \sym(\Omega(Y_0,\dots,Y_k)).$$
\end{definition}

	It is a straightforward consequence of the definition that $\group(N;N_0,\dots,N_k)=\smallgroup(N_0,\dots,N_k)$, and if $H\leq H'$ then $\group(H;N_0,\dots,N_k)\subset \group(H';N_0,\dots,N_k)$ (assuming that both sides are defined).

\begin{lemma}\label{why_group}
	$\group(H;N_0,\dots,N_k)$ is a group.
\end{lemma}

\begin{proof}
	Let $\sigma,\tau$ be permutations in $\sym(\Omega(Y_0,\dots,Y_k))$ satisfying conditions (a)-(d) in Definition~\ref{def:groups2}. We have to show that the permutation $\sigma\tau^{-1}$ also satisfies conditions (a)-(d). For conditions (a) and (b) this is obvious.

	For condition (c) let $\vec{n}\in \{0\}\times \omega^k$ be arbitrary. Then by Lemma~\ref{rho_product} we have $$\rho_{\vec{n}}(\sigma\tau^{-1})=\rho_{\psi(\tau^{-1})(\vec{\mathbf{n}})}(\sigma)\rho_{\vec{n}}(\tau^{-1})=\rho_{\psi(\tau^{-1})(\vec{\mathbf{n}})}(\sigma)(\rho_{\psi(\tau^{-1})(\vec{\mathbf{n}})}(\tau))^{-1}\in HH^{-1}=H.$$

	For condition (d) let $\vec{n},\vec{m}\in \{0\}\times \omega^k$. Then we have
\begin{align*}
\rho_{\vec{n}}(\sigma\tau^{-1})&= \rho_{\psi(\tau^{-1})(\vec{\mathbf{n}})}(\sigma)(\rho_{\psi(\tau^{-1})(\vec{\mathbf{n}})}(\tau))^{-1} \equiv \\
&\equiv\rho_{\psi(\tau^{-1})(\vec{\mathbf{m}})}(\sigma)(\rho_{\psi(\tau^{-1})(\vec{\mathbf{m}})}(\tau))^{-1}=\rho_{\vec{n}}(\sigma\tau^{-1}) \pmod{N}.
\end{align*}
\end{proof}

From this point on, we will use $\biggroup$ to denote $\aut(\Omega(Y_0,\dots,Y_k);Y_0\times\{0\},\bigeq,\smalleq)$. Continuing our analysis, we name a few sets:
\[
\begin{split}
T_i\coloneqq&\{\sigma\in  \biggroup\colon \exists n,m\in \omega (\psi(\sigma)=((n,i)(m,i)))  \wedge \forall j\forall n(\rho_n^j(\sigma)=\id(Y_j))\} \\
N_{n,i}\coloneqq&\{\sigma\in \biggroup\colon \psi(\sigma)=\id, \rho_n^i(\sigma)\in N_i, (\forall (m,j)\neq (n,i))(\rho_n^j(\sigma)=\id(Y_j))\},\\
H^*\coloneqq&\{\tau^*\colon \tau\in H\},
\end{split}
\]
where $\tau^*$ denotes the unique permutation $\sigma\in \biggroup$ which fixes every $\smalleq$ class setwise, and for which $\rho_{\vec{n}}(\sigma)=\tau$ holds for all $\vec{n}\in \{0\}\times \omega^k$.

	$T_i$ then contains those permutations in $\biggroup$ which transpose two copies of $Y_i$ and does not ``scramble'' any copies of any $Y_i$ (this is the condition $\rho_n^j(\sigma)=\id(Y_j)$). In particular, $\phi(\sigma)=\id(\{0,\ldots,k\})$. On the other hand, $N_{n,i}$ contains permutations that fix every $\smalleq$ class setwise and act as an element of $N_i$ on $Y_i\times\{n\}$ without moving anything else. 
 
\begin{lemma}\label{how_generated}
$\group(H;N_0,\dots,N_k)$ is locally generated by $\bigcup_{i=1}^kT_i\cup \bigcup_{i=0}^kN_{0,i}\cup H^*$.
\end{lemma}

\begin{proof}
	First we show that $\bigcup_{i=1}^kT_i\cup \bigcup_{i=0}^kN_{0,i}\cup H^*\subset \group(H;N_0,\dots,N_k)$. Let $\sigma\in \bigcup_{i=1}^kT_i\cup \bigcup_{i=0}^kN_{0,i}\cup H^*$. Then by definition $\sigma\in \biggroup$, and hence $\sigma$ automatically satisfies items (a) and (b) in Definition~\ref{def:groups2}. 
\begin{description}
\item[Case 1]{If $\sigma\in T_i$, then we have $\rho_{\vec{n}}=\id(Y)$ for all $\vec{n}\in \{0\}\times \omega^k$, and thus in this case items (c) and (d) are trivially satisfied.}
\item[Case 2]{If $\sigma\in N_{0,i}$, and $\vec{n}=(n_0,n_1,\dots,n_k)\in \{0\}\times \omega^k$ then we have $\rho_{\vec{n}}=\tau\cup \id(Y\setminus Y_i)\in N$ for some $\tau\in N_i$ if $n_i=0$, and $\rho_{\vec{n}}=\id(Y)$ otherwise. In any case, for all $\vec{n},\vec{m}\in \{0\}\times \omega^k$ we have $\rho_{\vec{n}}\in N\subset H$, so (c) is satisfied. Moreover, $\rho_{\vec{n}}(\sigma)(\rho_{\vec{m}}(\sigma))^{-1}\in \{\id\}\cup N\cup N^{-1}=N$, and so $\sigma$ satisfies (d).}
\item[Case 3]{Finally, if $\sigma=\tau^*$ then by definition for all $\vec{n},\vec{m}\in \{0\}\times \omega^k$ we have $\rho_{\vec{n}}(\sigma)$ is a uniquely determined element of $H$, so $\rho_{\vec{n}}(\sigma)(\rho_{\vec{m}}(\sigma))^{-1}=\id(Y)\in N$.}
\end{description}

	This proves $\bigcup_{i=1}^kT_i\cup \bigcup_{i=0}^kN_{0,i}\cup H^*$ is a subgroup of $\group(H;N_0,\dots,N_k)$. Now we show that in fact $\overline{\langle\bigcup_{i=1}^kT_i\cup \bigcup_{i=0}^kN_{0,i}\cup H^*\rangle}=\group(H;N_0,\dots,N_k)$.

\begin{claim}\label{claim:claim1}
The closure of the group $\langle\bigcup_{i=1}^kT_i\rangle$ contains every $\sigma\in \biggroup$ for which $\rho_n^j(\sigma)=\id(Y_j)$ for all $0\leq j\leq k$.
\end{claim}
\begin{proof}
Let $\perm:=\{\tau\in\sym(\{(0,0)\}\cup\{1,\ldots,k\}\times\omega)\colon \hat{\tau}\in \overline{\langle \bigcup_{i=1}^k T_i\rangle}\}$ where $\hat{\tau}$ is the unique element of $\biggroup$ such that for all $0\leq j\leq k$ we have $\rho_n^j(\hat{\tau})=\id(Y_j)$, and $\psi(\hat{\tau})=\tau$. We have to show that $\perm=\id(\{(0,0)\})\times \prod_{i=1}^k\sym(\{i\}\times\omega)$. The map $\tau\mapsto \hat{\tau}$ is a continuous homomorphism, so $\perm$ is a closed subgroup of $\id(\{(0,0)\})\times \prod_{i=1}^k\sym(\{i\}\times\omega)$. On the other hand, it follows directly from the definition of the sets $T_i$ that $\perm$ contains all the transpositions $(n,i)(m,i)$ for $1\leq i\leq k$ and $n,m\in \omega$. Since these transpositions locally generate the group $\id(\{(0,0)\})\times \prod_{i=1}^k\sym(\{i\}\times\omega)$ we obtain $\perm=\id(\{(0,0)\})\times \prod_{i=1}^k\sym(\{i\}\times\omega)$. 
\end{proof}

	Now we show that $\group(H;N_0,\dots,N_k)\subset \overline{\langle\bigcup_{i=1}^kT_i\cup \bigcup_{i=0}^kN_{0,i}\cup H^*\rangle}$. Let $\sigma\in \group(H;N_0,\dots,N_k)$ be arbitrary. Then by Claim~\ref{claim:claim1} we know that $\widehat{\psi(\sigma)}\in \overline{\langle\bigcup_{i=1}^kT_i\rangle}$. Therefore it is enough to show that $\sigma':=(\widehat{\psi(\sigma)})^{-1}\sigma\in \overline{\langle \bigcup_{i=0}^kN_{0,i}\cup H^*\rangle}$. By Claim~\ref{claim:claim1} we know that $(\widehat{\psi(\sigma)})^{-1}\in \overline{\langle\bigcup_{i=1}^kT_i\rangle}\subset \group(H;N_0,\dots,N_k)$. Therefore $\sigma'\in \group(H;N_0,\dots,N_k)$, and $\psi(\sigma')=\psi(\widehat{\psi(\sigma)})^{-1}\psi(\sigma)=\psi(\sigma)^{-1}\psi(\sigma)=\id$.

	Since $\sigma'\in \group(H;N_0,\dots,N_k)$ it follows that for all $\vec{n}\in \{0\}\times \omega^k$ we have $\rho_{\vec{n}}(\sigma')=(\rho_{\vec{0}}(\sigma'))^{-1}\rho_{\vec{n}}(\sigma)\in N$. Since $\sigma'$ also fixes every $\smalleq$-class setwise this means that $\sigma'$ can be written as the direct product of the groups $N_{n,i}\colon(n,i)\in \{(0,0)\}\cup \omega\times \{1,\dots,k\}$. In particular $\sigma'$ is locally generated by the groups $N_{n,i}$. On the other hand it is easy to verify that $N_{n,i}=\hat{\tau}^{-1}N_{0,i}\hat{\tau}$ where $\tau$ is the transposition $((0,i)(n,i))$. Thus $N_{n,i}\subset \langle N_{0,i}, T_i\rangle$. Therefore the direct product of the groups $N_{n,i}\colon (n,i)\in \{(0,0\}\cup \{1,\dots,k\}\times \omega$ is locally generated by $N_{0,i}\colon 0\leq i\leq k$ and $T_i\colon 1\leq i\leq k$, and hence $\sigma'\in \overline{\langle\bigcup_{i=1}^kT_i\cup \bigcup_{i=0}^kN_{i,0}\rangle}$ which finishes the proof of the lemma.
\end{proof}

	It is a well-known fact that every normal subgroup of $\sym(\omega)$ is either dense or trivial (see e.g.\ Proposition~1.4(2) in~\cite{MoonStalder} or Corollary~2.2 in~\cite{cameron1981normal}). This can be generalized to powers of $\sym(\omega)$ as follows.

\begin{lemma}\label{normal_power}
	Let $\perm$ be a normal subgroup of $\sym(\omega)^k$. Then either $\perm$ is dense in $\sym(\omega)^k$ or $\perm$ fixes every element of $\omega\times \{i\}$ for some $0\leq i\leq k-1$. 
\end{lemma}

\begin{proof}
	Let $\perm_i:=\perm|_{((\omega\times \{i\}))}$. Since $\perm$ is a normal subgroup of $\sym(\omega)^k$ it follows that $\perm_i$ is a normal subgroup of $\sym(\omega\times \{i\})$. Therefore $\perm_i$ is either dense in $\sym(\omega\times \{i\})$ or it is trivial for all $i<k$. If all groups $\perm_i$ are dense then 
\begin{align*}
\overline{\perm}\supset &\overline{\bigl\{\bigcup_{i<k}\alpha_i\colon \alpha_i\in \perm_i\bigr\}}\supset \\&\bigl\{\bigcup_{i<k}\alpha_i\colon \alpha_i\in \overline{\perm_i}\bigr\}=\bigl\{\bigcup_{i<k}\alpha_i\colon \alpha_i\in \sym(\omega\times \{i\})\bigr\}=\sym(\omega)^k,
\end{align*}
	that is $\perm$ is dense in $\sym(\omega)^k$. Let us suppose that this is not the case. Then $\perm_i$ is trivial for some $i$. Let $\perm_i':=\perm|_{\omega\times \{i\}}$. Then we have to show that $\perm_i'$ is also trivial.

	We show the contrapositive. Let us assume that $\perm_i'$ contains a nontrivial element. This means that there exists a permutation $\alpha\in \perm$ and $n\in \omega$ such that $\alpha((n,i))\neq (n,i)$. Let $m\in \omega$ such that $(m,i)$ is different from both $(n,i)$ and $\alpha((n,i))$. Let $\nu\in \sym(\omega\times k)$ denote the transposition flipping the elements $(n,i)$ and $(m,i)$. Then $\nu\in \sym(\omega)^k$ and since $\perm\triangleleft \sym(\omega)^{k}$ we have $\beta:=\nu^{-1}\alpha^{-1}\nu\alpha\in \perm$. On the other hand $$\beta((n,i))=\nu^{-1}\alpha^{-1}\nu\alpha((n,i))=\nu\alpha^{-1}(\alpha((n,i)))=\nu((n,i))=(m,i)\neq (n,i),$$ and since $\nu$ fixes every element outside $\omega\times \{i\}$ so does $\beta$. This means that $\beta\in \perm_{\omega\times (k\setminus \{i\})}$ and thus $\beta|_{\omega\times \{i\}}\in \perm_i$ which contradicts our assumption that $\perm_i$ is trivial.
\end{proof}

\begin{lemma}\label{normal_subgroups}
	The group $\smallgroup(N_0,\dots,N_k)$ has no nontrivial oligomorphic closed normal subgroup.
\end{lemma}

\begin{proof}
	Let $G$ be an oligomorphic closed normal subgroup of $\smallgroup(N_0,\dots,N_k)$. We have to show that $G=\smallgroup(N_0,\dots,N_k)$.

	In this proof we use the notation from the proof of Lemma~\ref{how_generated}. 

	Since $N_0=\id(Y_0)$ it follows that $\smallgroup(N_0,\dots,N_k)$, and hence also $G$ fixes $Y_0\times \{0\}$ pointwise.

	Let $\perm:=\{\psi(\sigma)|_{\omega \times \{1,\dots,k\}}\colon \sigma\in G\}$. Since $G$ preserves the relations $\bigeq$ and $\smalleq$, and $\psi$ is a homomorphism it follows that $\perm$ is a subgroup of $\sym(\omega)^{[k]}$ where $[k]=\{1,\dots,k\}$. We claim that $\perm$ is in fact a normal subgroup of $\sym(\omega)^{[k]}$. Indeed if $\sigma\in G$, and $\tau\in \sym(\omega)^{[k]}$, then $\hat{\tau}\in \smallgroup(N_0,\dots,N_k)$ according to Lemma~\ref{how_generated}, and then since $G$ is normalized by $\smallgroup(N_0,\dots,N_k)$ it follows that $\tau^{-1}\sigma\tau=\psi_i(\hat{\tau}^{-1}\sigma\hat{\tau})\in \perm$. Therefore by Lemma~\ref{normal_power} we know that either $\perm$ is a dense subgroup of $\sym(\omega)^{[k]}$ or there exists some $i\in \{1,\dots,k\}$ such that $\perm$ fixes every element in $\omega\times \{i\}$. Since $\psi$ is a homomorphism the latter condition would imply that $G$ fixes the all the sets $Y_i\times \{n\}$ setwise. In particular $G$ has infinitely many orbits, a contradicition. Therefore $\perm$ is dense in $\sym(\omega)^{[k]}$. $(*)$.

	Let us now fix a natural number $n\geq 3$, and let $i\in \{1,\dots,k\}$. Then there exists a $\sigma\in G$ such that $\psi(\sigma)|_{\{0,\dots,n\}\times \{1,\dots,k\}}=((0,i),(n,i))$. Let $\mu\in \aut(\Omega(Y_0,\dots,Y_k);Y_0,\bigeq,\smalleq)$ such that $\psi(\mu)=(0,i),(n,i))$ and $\rho_n^j(\mu)=\id(Y_j)$ for all $j$. Then $\mu\in T_i$. For a permutation $\tau\in N_i$ we denote by $\tilde{\tau}$ the permutation of $\Omega(Y_0,\dots,Y_k)$ which maps $(a,0)$ to $(\tau(a),0)$ for $a\in Y_i$ and fixes every other element. Then by definition $\tilde{\tau}\in N_{0,i}$. By Lemma~\ref{how_generated} we know that $T_i,N_{i,0}\subset \smallgroup(N_0,\dots,N_k)$. Thus $\smallgroup(N_0,\dots,N_k)$ contains $\mu$ and $\tilde{\tau}$ for all $\tau\in N_i$. Now since $G$ is normalized by $\smallgroup(N_0,\dots,N_k)$ it follows that $G\ni \nu_{\tau}:=\sigma^{-1}\tilde{\tau}^{-1}\sigma\tilde{\tau}$. Then we have $$\psi(\nu_{\tau})=\psi(\sigma)^{-1}\psi(\tilde{\tau})^{-1}\psi(\sigma)\psi(\tilde{\tau})=\psi(\sigma)^{-1}\psi(\sigma)=\id.$$

	On the other hand by using Lemma~\ref{rho_product} we obtain
\begin{align*}
\rho_0^i(\nu_{\tau})=&(\rho_{\psi_i(\sigma^{-1}\tilde{\tau}^{-1}\sigma\tilde{\tau})(0)}^i(\sigma))^{-1}(\rho_{\psi_i(\tilde{\tau}^{-1}\sigma\tilde{\tau})(0)}^i(\tilde{\tau}))^{-1}\rho_{\psi_i(\tilde{\tau})(0)}^i(\sigma)\rho_0^i(\tilde{\tau})=\\=&(\rho_0^i(\sigma))^{-1}(\rho_n^i(\tilde{\tau}))^{-1}\rho_n^i(\sigma)\rho_0^i(\tilde{\tau})=(\rho_0^i(\sigma))^{-1}\rho_0^i(\sigma)\tau=\tau,
\end{align*}

	and $\psi(\nu_{\tau})$ fixes every element of $\{1,\dots,n-1\}\times \{1,\dots,k\}$. This implies that $\nu_{\tau}$ and $\tilde{\tau}$ agrees on $Y_0\cup \bigcup_{i=1}^k(Y_i\times \{0,\dots,n-1\})$.

	Since this is true for all $n\geq 3$, and the group $G$ is closed it follows that $\tilde{\tau}\in G$. Therefore $N_{0,i}\leq G$. It follows similarly that $N_{n,i}\leq G$ for all $(n,i)\in \omega\times \{1,\dots,k\}$. Since $G$ is closed this implies that $G$ also contains the direct product of the groups $N_{n,i}$ $(\dag)$. Thus by using Lemma~\ref{how_generated} it remains to show that $T_i\subset G$ for all $1\leq i\leq k$.

	We claim that for all $\sigma\in G$ the permutation $\widehat{\psi(\sigma)}$ is also contained in $G$. Indeed let us observe that $$\psi(\widehat{\psi(\sigma)}\sigma^{-1})=\psi(\widehat{\psi(\sigma)})(\psi(\sigma))^{-1}=\psi(\sigma)(\psi(\sigma)^{-1})=\id.$$ Since $\widehat{\psi(\sigma)}\sigma^{-1}\in \smallgroup(N_0,\dots,N_k)$ this means that the permutation $\widehat{\psi(\sigma)}\sigma^{-1}$ is contained in the direct product of the groups $N_{n,i}$. Then by $(\dag)$ we obtain that $\widehat{\psi(\sigma)}\sigma^{-1}\in G$, and thus $\widehat{\psi(\sigma)}\in G$. Combining this with $(*)$ we obtain that $\{\tau\in \sym(\omega)^{[k]}\colon \hat{\tau}\in G\}$ is a dense subgroup of $\sym(\omega)^{[k]}$. On the other hand since the map $\tau\mapsto \hat{\tau}$ is continuous it follows that $\{\tau\in \sym(\omega)^k\colon \hat{\tau}\in G\}$ is also closed. Hence we obtain that in fact $\{\tau\in \sym(\omega)^{[k]}\colon \hat{\tau}\in G\}=\sym(\omega)^{[k]}$, that is $\hat{\tau}\in G$ for all $\tau\in \sym(\omega)^{[k]}$. This implies in particular that $T_i\subset G$ for all $1\leq i\leq k$.
\end{proof}

\begin{lemma}\label{small_normal_subgroup}
	The group $\smallgroup(N_0,\dots,N_k)$ is a normal subgroup of $\group(H;N_0,\dots,N_k)$, and $\group(H;N_0,\dots,N_k)/\smallgroup(N_0,\dots,N_k)\simeq H/N$.
\end{lemma}

\begin{proof}
	Condition (d) in Definition~\ref{def:groups2} says that for a permutation $\sigma\in \group(H;N_0,\dots,N_k)$ all the permutations of the form $\rho_{\vec{n}}(\sigma)$ are the same modulo $N$. Then we define $h(\sigma)$ to be the coset of $N$ containing all these elements. Then $h$ defines a homomorphism from $\group(H;N_0,\dots,N_k)$ to $H/N$, and the kernel of this homomorphism is $\smallgroup(N_0,\dots,N_k)$. This implies the statement of the lemma.
\end{proof}

\begin{lemma}\label{modulo_equal}
	Let $G\subset\sym(\Omega(Y_0,\ldots,Y_k))$ be a closed permutation group which preserves the set $Y_0\times\{0\}$ and the equivalence relations $\smalleq$ and $\bigeq$. Let us assume moreover that $\smallgroup(N_0,\dots,N_k)\leq G$, and for all $i=1,\dots,k$ we have $G_{((Y_i\times \{0\}))}=N_i\times \{\id(\{0\})\}$. Then for all $\vec{n},\vec{m}\in \{0\}\times \omega^k$ and $\sigma\in G$ the permutation $\rho_{\vec{n}}(\sigma)$ and $\rho_{\vec{m}}(\sigma)$ are equal modulo $N$.
\end{lemma}

\begin{proof}
	It is enough to show the statement of the lemma for $\vec{m}=\vec{0}=(0,\dots,0)$ and any $\vec{n}=(n_0,\dots,n_k)$. By definition both $\rho_{\vec{0}}$ and $\rho_{\vec{n}}$ maps $Y_i$ to $Y_{\phi(\sigma)(i)}$. Therefore $\tau\coloneqq\rho_{\vec{n}}(\sigma)^{-1}\rho_{\vec{0}}(\sigma)$ fixes the set $Y_0,\dots,Y_k$ setwise.

	Now let $\tau_i\coloneqq\tau|_{Y_i}=\rho_{\vec{0}}(\sigma)^{-1}\rho_{\vec{n_i}}(\sigma)$. We have to show that $\tau_i\in N_i$. If $n_i=0$ then there is nothing to prove. Otherwise let $\nu\in \sym(\bigcup_{i=0}^n)$ which swaps $(0,i)$ and $(n,i)$ and fixes everything else. Then using the notation of Lemma~\ref{how_generated} we have $\hat{\nu}\in T_i\subset \smallgroup(N_0,\dots,N_k)\subset G$. Now let $\mu:=\hat{\nu}^{-1}\sigma^{-1}\hat{\nu}\sigma$. Then $\mu$ fixes every element outside $Y_i\times \{0,n\}$ and by Lemma~\ref{rho_product} we have
\begin{align*}
\rho_0^i(\mu)=&\rho_0^i(\hat{\nu}^{-1}\sigma^{-1}\hat{\nu}\sigma)=\\
&(\rho_{\psi_i(\hat{\nu}^{-1}\sigma^{-1}\hat{\nu}\sigma)(0)}^i(\hat{\nu}))^{-1}(\rho_{\psi_i(\sigma^{-1}\hat{\nu}\sigma)(0)}^i(\sigma))^{-1}\rho_{\psi_i(\sigma)(0)}^i(\hat{\nu})\rho_0^i(\sigma)=\\
&(\rho_n^i(\sigma))^{-1}\rho_0^i(\sigma)=\tau_i.
\end{align*}
	
	Now by conjugating $\sigma'$ by the permutation $\widehat{\nu_{n'}}$ where $\nu_{n'}=((i,n)(i,n'))$ we obtain that for all $n'\in \omega$ there exists a permutation $\mu'_{n'}$ in $G$ such that $\rho_0^i(\mu_{n'}')=\tau_i$, and $\mu'_{n'}$ fixes every element outside $Y_i\times \{0,n'\}$. Since $G$ is closed this implies that $G$ contains a permutation $\mu'$ such that $\rho_0^i(\mu')=\tau$ and $\mu'$ fixes every element outside $Y_i\times \{0\}$. Then the permutation $(\tau_i,\id(\{0\}))$ is contained in $G_{(((Y_i\times \{0\}))}=N_i\times \{\id(\{0\})\}$, and thus $\tau_i\in N_i$, and this is what we wanted to show.
\end{proof}

\begin{lemma}\label{n_normal_subgroup}
	Let $G$ be as in Lemma~\ref{modulo_equal}, and let $H:=\{\rho_{\vec{0}}(\sigma)\colon \sigma\in G\}$. Then $N\triangleleft H$.
\end{lemma}

\begin{proof}
	Let $\sigma\in N$, $\tau\in H$. We have to show that $\tau^{-1}\sigma\tau\in N$. Let $\sigma_0$ be the permutation of $\Omega(Y_0,\dots,Y_k)$ which maps $(a,0)$ to $(\sigma(a),0)$ for all $a\in Y$, and fixes every other element. Then by definition $\sigma_0\in \smallgroup(N_0,\dots,N_k)\subset G$. Let $\tilde{\tau}\in G$ so that $\rho_{\vec{0}}(\tilde{\tau})=\tau$. Then by permuting the $\smalleq$ classes of $\Omega(Y_0,\dots,Y_k)$ with elements in $T_i$ (see Lemma~\ref{how_generated}) we can assume without loss of generality that $\psi(\tilde{\tau})$ fixes all $(0,i)\colon i=0,\dots,k$. In this case the permutation $\sigma_0':=\tilde{\tau}^{-1}\sigma_0\tilde{\tau}\in G$ fixes every element outside $\bigcup_{i=0}^kY_i\times \{0\}$, and by Lemma~\ref{rho_product} we have
\[\rho_{\vec{0}}(\sigma_0')=\rho_{\vec{0}}(\tilde{\tau}^{-1}\sigma_0\tilde{\tau})=(\rho_0(\tilde{\tau}))^{-1}\rho_0(\sigma_0)\rho_0(\tilde{\tau})=\tau^{-1}\sigma\tau.
\]
	This implies that $\tau^{-1}\sigma\tau\in N$.
\end{proof}

\begin{lemma}\label{desc}
	Let $G_1$ and $G_2$ be as in Lemma~\ref{modulo_equal}, and let $H_i:=\{\rho_{\vec{0}}(\sigma)\colon \sigma\in G\}$. If $H_1=H_2$ then $G_1=G_2$.
\end{lemma}

\begin{proof}
	It follows from our conditions that the groups $G_1$ and $G_2$ both contain the group $\smallgroup(N_0,\dots,N_k)=\id(Y_0\times \{0\})\times \prod_{i=1}^k(N_i\wr \sym(\omega))$.

	Lemma~\ref{modulo_equal} implies that for all $\sigma\in G_1$ the coset $\rho_{\vec{n}}(\sigma)N$ does not depend on $\vec{n}$. Let us denote this coset by $h(\sigma)$. By Lemma~\ref{n_normal_subgroup} we know that $N$ is a normal subgroup of $H_1$. Then it follows easily from Lemma~\ref{rho_product} that $h$ defines a homomorphism from $G_1$ to $H_1/N$, and the kernel of $h$ is exactly $\smallgroup(N_0,\dots,N_k)$.

	Using this if $\sigma\in G_1$ then by definition it follows that $\rho_{\vec{0}}(\sigma')=\rho_{\vec{0}}(\sigma)$ for some $\sigma'\in G_2$. This implies $h(\sigma)=h(\sigma')$, and thus $h(\sigma(\sigma')^{-1})=h(\sigma)(h(\sigma'))^{-1}=\id$, that is $\rho_{\vec{n}}(\sigma(\sigma')^{-1})$ for all $\vec{n}\in 0\times \omega^k$. This means exactly that $\sigma(\sigma')^{-1}\in \smallgroup(N_0,\dots,N_k)$. Therefore $\sigma=(\sigma(\sigma')^{-1})\sigma'\in \smallgroup(N_0,\dots,N_k)G_2=G_2$. Hence $G_1\subset G_2$. The other containment follows analogously.
\end{proof}

\begin{lemma}\label{recursive_description}
	Let $G$ be as in Lemma~\ref{modulo_equal}. Then there exists a permutation group $H$ acting on $Y$ which normalizes $N$, and such that $G=\group(H;N_0,\dots,N_k)$.
\end{lemma}

\begin{proof}
	Let $H=H_1:=\{\rho_{\vec{0}}(\sigma)\colon \sigma\in G\}$. By Lemma~\ref{desc} we already know that $H$ normalizes $N$. Let $G_1=G, G_2=\group(H;N_0,\dots,N_k)$, and $H_2=\{\rho_{\vec{0}}(\sigma)\colon \sigma\in G_2\}$. Then $H_2=H=H_1$, and thus Lemma~\ref{desc} implies $G=G_1=G_2=\group(H;N_0,\dots,N_k)$.
\end{proof}

\begin{definition}
	We define the classes of permutation groups $\gone_{-1}',\gone_0',\gone_1',\dots$ recursively as follows.
\begin{itemize}
\item $\gone_{-1}'=\emptyset$,
\item For $n\geq 0$ we define $G\in \gone_n'$ if and only if there exists $N_0,\dots,N_k,N,H$ as in Definition \ref{def:groups2} such that $N_1,\dots,N_k,H\in \gone_{n-1}'$ and $G$ is isomorphic to $\group(H;N_0,\dots,N_k)$.
\end{itemize}
	We define $\gone'$ as the union of all classes $\gone_n'\colon n\in \omega$.
\end{definition}

\begin{remark}
	Later we will see that in fact the classes $\gone_n'$ are closed under finite direct products. Therefore we could drop the assumption that $N\in \gone_{n-1}'$ in the recursion. This fact can also be derived from the definition directly. We will also see that we can drop the condition $H\in \gone_{n-1}'$ as well.
\end{remark}

\begin{lemma}\label{gone_oligo}
	Every group in $\gone'$ is oligomorphic.
\end{lemma}

\begin{proof}
	We show that every group in $\gone_n'$ is oligomorphic by induction on $n$. For $n=-1$ there is nothing to show. For the induction step let us assume that $G=\group(H;N_0,\dots,N_k)$ for some $H,N_0,\dots,N_k\in \gone_{n-1}'$. Then $G\geq \smallgroup(N_0,\dots,N_k)$. By item (ii) of Lemma~\ref{wr_closed} we know that the groups $N_i\wr \sym(\omega)\colon 1\leq i \leq k$ are oligomorphic. By the discussion in Section~\ref{sect:direct} we know that finite direct products of oligomorphic groups is oligomorphic. Therefore the group $\smallgroup(N_0,\dots,N_k)=\id(Y_0\times \{0\})\times \prod_{i=1}^k(N_i\wr \sym(\omega))$ is oligomorphic. Hence $G$ is also oligomorphic.
\end{proof}

\begin{lemma}\label{finite_index_cl}
	Let $G\in \gone'$. Then $G$ has finitely many closed normal subgroups contained in $\gone'$, and all of these have a finite index in $G$.
\end{lemma}

\begin{proof}
	We show the statement of the lemma for $\gone_n'$ by induction on $n$. For $n=-1$ there is nothing to show. For the induction step let us assume that $G=\group(H;N_0,\dots,N_k)$ for some $H,N_0,\dots,N_k,N=\prod_{i=0}N_i\in \gone_{n-1}'$. By Lemma~\ref{small_normal_subgroup} it follows that $|G:\smallgroup(N_0,\dots,N_k)|=|H:N|$. By definition $N$ is a normal subgroup of $H\in \gone_{n-1}'$. Since $N,H\in \gone_{n-1}$ this implies that $|H:N|$ is finite by the induction hypothesis. Therefore $|G:\smallgroup(N_0,\dots,N_k)|$ is finite. Now let $G_1\in \gone'$ be a closed normal subgroup of $G$. By Lemma~\ref{gone_oligo} we know that $G_1$ is oligomorphic. Let $G_2:=G_1\cap \smallgroup(N_0,\dots,N_k)$. Then by the second isomorphism theorem we obtain 
\begin{multline*}
G_1/G_2=G_1/(G_1\cap \smallgroup(N_0,\dots,N_k))\simeq G_1(N_0,\dots,N_k)/\smallgroup(N_0,\dots,N_k)\leq \\ \leq G/\smallgroup(N_0,\dots,N_k).
\end{multline*}
	In particular $G_1/G_2$ is finite. We know that a finite index subgroup of an oligomorphic group is oligomorphic (see for instance Proposition 6.23 in~\cite{bodirsky2018structures}). Therefore $G_2$ is also oligomorphic. $G_2$ is also closed since it is an intersection of two closed groups. Thus by Lemma~\ref{normal_subgroups} it follows that $\smallgroup(N_0,\dots,N_k)\leq G_2\leq G$. Since $|G:\smallgroup(N_0,\dots,N_k)|< \infty$ this implies both statements of the lemma.
\end{proof}

\begin{lemma}\label{everything_closed}
	Every group in $\gone'$ is closed.
\end{lemma}

\begin{proof}
	We show by induction on $n$ that every group in $\gone_n'$ is closed. For $n=-1$ there is nothing to prove. For the induction step let us assume that $G=\group(H;N_0,\dots,N_k)$ for some $H,N_0,\dots,N_k,N=\prod_{i=0}N_i\in \gone_{n-1}'$. By item (i) of Lemma~\ref{wr_closed} we know that the groups $N_i\wr \sym(\omega)\colon 1\leq i \leq k$ are closed. Thus the direct product $\smallgroup(N_0,\dots,N_k)=\id(Y_0\times \{0\})\times \prod_{i=1}^k(N_i\wr \sym(\omega))$ is also closed (see Section~\ref{sect:direct}). Since $\smallgroup(N_0,\dots,N_k)\triangleleft \group(H;N_0,\dots,N_k)$ it follows that $|G:\smallgroup(N_0,\dots,N_k)|$ is finite. This means that $G$ is a finite union of cosets of $\smallgroup(N_0,\dots,N_k)$, and thus a finite union of closed sets. Hence $G$ is closed.
\end{proof}

	By Lemma~\ref{everything_closed} we can remove the condition ``closed'' from the formulation of Lemma~\ref{finite_index_cl}.

\begin{corollary}\label{finite_index}
	Let $G\in \gone'$. Then $G$ has finitely many normal subgroups contained in $\gone'$, and all of these have a finite index in $G$.
\end{corollary}

\begin{lemma}\label{easy_containment}
	$\gone_n'\subset \gone_n$ for all $n\in \omega$.
\end{lemma}

\begin{proof}
	We prove the lemma by induction on $n$. For $n=-1$ the statement is trivial. For the induction step let us assume that $G=\group(H;N_0,\dots,N_k)$ for some $H,N_0,\dots,N_k,N=\prod_{i=0}N_i\in \gone_{n-1}'$. By Lemmas~\ref{inf_copies2} and~\ref{direct_gone} it follows that $\prod_{i=1}^k(N_i\wr \sym(\omega))\times N_0=\smallgroup(N_0,\dots,N_k)\in \gone_n$. Since $\gone_n$ is closed under taking closed supergroups it follows that $\overline{G}\in \gone_n$. Finally we have $\overline{G}=G$ by Lemma~\ref{everything_closed}.
\end{proof}

	Now we show that in fact $\gone_n=\gone_n'$ for all $n$. We prove this by induction on $n$. For the induction step we first show that if $\gone_{n-1}=\gone_{n-1}'$, and if $G\in \gone$ with $\omega$-partition $\partition=(Y_0, \bigeq, \smalleq)$ with components $G_0,\dots,G_k$ then $G$ has a subgroup $S$ which is isomorphic to a group of the form $\id(Y_0)\times \prod_{i=1}^m(G_i\wr \sym(\omega))$. First we show this in the case when $K=\emptyset$, and $\bigeq$ is the universal congruence. The proof of this step is based on the argument presented in Case 2 of the proof of Lemma 3.3, part (i) in~\cite{lachlan1992}.

\begin{lemma}\label{splitting}
	Let us assume that $\gone_{n-1}=\gone_{n-1}'$.

	Let $G\in \gone_n$ be a permutation group with $\dom(G)=X$, and let $\partition=(\emptyset, X^2, \smalleq)$ be the $\omega$-partition of $G$. Let $Y^0\in X/\smalleq$. Then there exists a map $p\colon X\rightarrow Y^0$ such that 
\begin{enumerate}
\item $p|_{Y^0}=\id(Y^0)$,
\item $p$ defines a bijection from $Y$ to $Y_0$ for all $Y\in X/\smalleq$, and
\item for all $\sigma\in \sym(\omega)$ there exists a permutation $\hat{\sigma}\in G$ which acts on $X/\smalleq$ as $\sigma$, and such that $p(\hat{\sigma}(u))=p(u)$ for all $u\in X$.
\end{enumerate}
\end{lemma}

\begin{proof}
	Let $Y_0,Y_1,\dots$ be an enumeration of all $\smalleq$-classes of $X$. For a subset $S$ of $\omega$ we define $Y_S:=\bigcup_{j\in S}(Y_j)$. Let 
$$N_S:=(G_{X\setminus Y_S} \cap \bigcap_{j\in S}G_{\{Y_j\}})|_{Y_S}, \text{ and } H_S:=(\bigcap_{j\in S}G_{\{Y_j\}})|_{Y_S}.$$ We use the notation $N_j,H_j$ for the groups $N_{\{j\}}=G_{((Y^j))}$ and $H_{\{j\}}=G_{(Y^j)}$, respectively. It follows from the definition that for all $S\subset \omega$ the group $N_S$ is closed, and $N_S$ is a normal subgroup of $H_S$. Also we have $\prod_{j\in S}N_j\subset N_S$. By definition every $N_j$ is a component of the $\omega$-partition $\partition$. Therefore $N_j\in \gone_{n-1}$ for all $j\in \omega$. Then Lemma~\ref{direct_gone} implies that if $S$ is finite then $\prod_{j\in S}N_j\in \gone_{n-1}$. This also implies that the groups $N_S$ and $\overline{H_S}$ are contained $\gone_{n-1}$ since $\gone_{n-1}$ is closed under taking closed supergroups. Since $N_S$ is closed, and $N_S\triangleleft H_S$ it follows that $h^{-1}N_Sh\in N_S$ for all $h\in \overline{H_S}$. Therefore $N_S$ is also a normal subgroup of $\overline{H_S}$. Since $\gone_{n-1}=\gone_{n-1}'$ this implies that $|\overline{H_S}:N_S|$ is finite by Corollary~\ref{finite_index}. Therefore $|H_S:N_S|$ is also finite (and thus in fact $H_S$ is closed). This also implies that the group $G_{X\setminus Y_S} \cap \bigcap_{j\in S}G_{\{Y_j\}}$ is a finite index subgroup of $\bigcap_{j\in S}G_{\{Y_j\}}$. The group $G$ acts on $X/\smalleq$ as the full symmetric group for all $i$. This means that $\bigcap_{j\in S}G_{\{Y_j\}}$ acts on $(X\setminus Y_S)/\smalleq$ as the full symmetric group. Since $\sym(\omega)$ has nontrivial finite index subgroup it follows that the group $G_{X\setminus Y_S} \cap \bigcap_{j\in S}G_{\{Y_j\}}$ acts on $(X\setminus Y_S)/\smalleq$ as the full symmetric group as well. In particular for any $j_1,j_2\in \omega\setminus S$ there exists a permutation $\mu\in G$ such that
\begin{itemize}
\item $\mu$ fixes every $Y_j\colon j\neq j_1,j_2$ setwise
\item $\mu$ switches $Y_{j_1}$ and $Y_{j_2}$, and
\item $\mu|_{Y_S}\in N_S$.
\end{itemize}

	Let $K_l:=(G_{\{Y_0\}}\cap G_{Y_{\{1,\dots,l\}}})|_{X_0}$. Then $K_0=H_0\supset K_1\supset \dots$, and since $G$ is closed it follows that $\bigcap_{l=0}^\infty{K_l}=N_0$. Since $|H_0:N_0|< \infty$ it follows that there exists an index $l$ such that $K_l=N_0$, and hence $K_{l'}=N_0$ for all $l'\geq l$. For all subsets $S\subset \omega\setminus \{0\}$ with $|S|\geq l$ there exists a permutation $\nu$ in $G$ which fixes $Y_0$, and which maps $Y_{\{1,\dots,|S|\}}$ into $Y_S$. Then $\nu$ maps $G_{\{X_0\}}\cap G_{Y_{\{1,\dots,|S|\}}}$ to $G_{\{X_0\}}\cap G_{Y_S}$, and thus it maps $N_0=K_l$ to $K_S:=(G_{\{X_0\}}\cap G_{Y_S})|_{Y_0}$. On the other hand $\nu$ fixes $N_0=G|_{((Y_0))}$ since $\nu$ fixes $Y_0$. Therefore $K_S=N_S$ for all $S \subset \omega\setminus \{0\}$ with $|S|\geq l$. $(*)$

	Let $j\in \omega\setminus 0$. We claim that for every finite set $S\subset \omega\setminus \{0,j\}$ there exists a permutation $\sigma_S$ of $G$ such that 
\begin{itemize}
\item $\sigma_S$ fixes $Y_S$ pointwise,
\item $\sigma_S$ switches $Y_0$ and $Y_j$, and 
\item $\sigma_S^2|_{Y_0\cup Y_j}=\id(Y_0\cup Y_j)$.
\end{itemize}
	We can assume without loss of generality that $|S|\geq l$. Then we have seen that there exists a permutation $\mu_S\in G_{S}$ that switches $Y_0$ and $Y_j$. By $(*)$ it follows that $(\mu_S)^2|_{Y_0}\in N_0$. This means that there exists a $\tau\in G|_{X\setminus Y_0}$ such that $\tau|_{Y_0}=\sigma^2|_{Y_0}$. Let $\sigma_S:=\tau^{-1}\mu_S$. Then $\sigma_S$ also switches $Y_0$ and $Y_j$, and it also fixes $Y_S$ pointwise. Moreover we have
\begin{align*}
(\sigma_S^2)|_{Y_0}=&(\tau^{-1}\mu_S)^2|_{Y_0}=(\tau^{-1}\mu_S\tau^{-1}|_{\mu_S(Y_0)}\mu_S)|_{Y_0}=\\&(\tau^{-1}\mu_S^2)|_{Y_0}=\tau|_{Y_0}^{-1}(\mu_S)^2|_{Y_0}=\id(Y_0).
\end{align*}
	It follows similarly that $(\tau^{-1}\sigma_S^2)|_{Y_j}=\id(Y_j)$. This means that $\sigma_S^2|_{Y_0\cup Y_j}=\id(Y_0\cup Y_j)$.

	Now let $S_0$ be such any finite subset with $|S_0|\geq l$. %We define $p|_{Y_j}:=\sigma_{S_0}|_{Y_j}$. 
We claim that for every finite subset $S$ of $\omega\setminus \{0,j\}$ there exists a map $\sigma'_S\in G|_{X^{S'}}$ such that $\sigma_S'|_{Y_0\cup Y_j}=\sigma_{S_0}|_{Y_0\cup Y_j}$. ($\dag$) This is enough to show in the case when $S'\supset S_0$. In this case $\sigma_{S_0}\sigma_S^{-1}$ fixes $Y_0$ setwise, and $Y_{S_0}$ pointwise. Therefore by $(*)$ we know that there exists a $\tau'\in G|_{X\setminus Y_0}$ such that $\tau'|_{Y^0}=(\sigma_{S_0}\sigma_S^{-1})|_{Y_0}$. By symmetrical argument we can also find a $\tau''\in G|_{X\setminus Y^j}$ such that $\tau''|_{Y_j}=(\sigma_{S_0}\sigma_S^{-1})|_{Y_j}$. Let $\tau:=\tau'\tau''$. Then $\tau|_{Y_0\cup Y_j}=(\sigma_{S_0}\sigma_S^{-1})|_{Y_0\cup Y_j}$. Let $\sigma_S':=\tau\sigma_S$. Then $\sigma_S'$ switches $Y_0$ and $Y_1$, fixes $S$ pointwise, and $$(\sigma_S')|_{Y_0\cup Y_j}=\tau_{Y_0\cup Y_j}\sigma_S|_{Y_0\cup Y_j}=(\sigma_{S_0}\sigma_S^{-1})|_{Y_0\cup Y_j}\sigma_S)|_{Y_0\cup Y_j}=\sigma_{S_0}|_{Y_0\cup Y_j}.$$ 

	Since $G$ is closed it also follows from $(\dag)$ that there exists a $\sigma\in G$ such that $\sigma|_{Y_0\cup Y_j}=\sigma_S|_{Y_0\cup Y_j}$ and $\sigma$ fixes every element outside $Y_0\cup Y_j$. In particular $\sigma$ switches the classes $Y_0$ and $Y_j$, and $\sigma^2|_{Y_0\cup Y_j}=\id(Y_0\cup Y_j)$, and thus $\sigma^2=\id(X)$. 
	
	We define $p$ on $Y_j$ as $\sigma|_{Y_j}$ for the map $\sigma$ as above defined for $j$. We extend $p$ to $Y_0$ as the identity map. Then $p$ satisfies conditions (1) and (2) of the lemma. We show that it also satisfies condition (3). 

	For a permutation $\sigma\in \sym(\omega)$ there exists a unique permutation $\hat{\sigma}\in \sym(X)$ such that $p(\hat{\sigma}(u))=p(u)$ for all $u\in X$, and $\hat{\sigma}/\smalleq=\sigma$. We have to show that $\hat{\sigma}\in G$ for all $\sigma\in \sym(\omega)$. By the definition $p$ it follows that $\hat{\sigma}\in G$ for every transposition $\sigma=(0j)\in \sym(\omega)$. The map $\sigma\mapsto \hat{\sigma}$ is a continuous group homomorphism. Therefore $\{\sigma\in \sym(\omega)\colon \hat{\sigma}\in G\}$ is a closed subgroup of $\sym(\omega)$. Since the transpositions of the from $(0j)\colon j\in \omega\setminus \{0\}$ locally generate $\sym(\omega)$ this implies that in fact $\{\sigma\in \sym(\omega)\colon \hat{\sigma}\in G\}=\sym(\omega)$ i.e., $\hat{\sigma}\in G$ for all $\sigma\in \sym(\omega)$. This finishes the proof of the lemma. 
\end{proof}

\begin{corollary}\label{wreath_product}
	Let us assume that $\gone_{n-1}=\gone_{n-1}'$.

	Let $G\in \gone_n$ be a permutation group with $\dom(G)=X$, and let $\partition=(\emptyset, X^2, \smalleq)$ be the $\omega$-partition of $G$. Let $Y^0\in X/\smalleq$. Then $G$ contains a group $G^*$ so that $\partition$ is also an $\omega$-partition of $G^*$ with the same components, and $G^*$ is isomorphic to $G_{((Y^0))}\wr \sym(\omega)$.
\end{corollary}

\begin{proof}
	Let $p$ be as in the statement of Lemma~\ref{splitting}, and let $Y_0,Y_1,\dots$ be the $\smalleq$-classes in $X$. We define $e\colon X\rightarrow Y_0\times \omega, u\mapsto (p(u),q(u))$ where $q(u)$ denotes the unique $j$ such that $u\in Y_j$. Then $e$ is a bijection between $X$ and $Y_0\times \omega$. Let $N_i:=G_{((Y^i))}$. Then since $G$ is closed it is clear $G$ contains the direct $\prod_{i=0}^\infty N_i$. The image of this direct product according to $\iota$ is exactly $N_0^\omega$. Using the map $\sigma\mapsto \hat{\sigma}$ as in the proof of Lemma~\ref{splitting} it is easy to see that the maps $\sigma\mapsto e(\hat{\sigma})$ maps $\sym(\omega)$ into $\id(Y_0)\wr \sym(\omega)$. Therefore the image of $e$ contains $Y_0\wr \sym(\omega)$. We put $G^*:=e^{-1}(N_0\wr \sym(\omega))$. Then $G^*$ satisfies the conditions of the lemma.
\end{proof}

\begin{lemma}\label{wreath_product_gen}
	Let us assume that $\gone_{n-1}=\gone_{n-1}'$.

	Let $G\in \gone_n$ be a permutation group with $\dom(G)=X$, and let $\partition=(K, \bigeq, \smalleq)$ be the $\omega$-partition of $G$. Let $X_1,\dots,X_k$ denote the $\bigeq$-classes. Let $Y_i\in X_i/\smalleq$. Let $N_i:=G|_{((Y_i))}$. Then $G$ contains a group $G^*$ so that $\partition$ is also an $\omega$-partition of $G_0$ with the same components as $G$, and $G_0$ is isomorphic to $\smallgroup(\id(Y_0),N_1,\dots,N_k)$.
\end{lemma}

\begin{proof}
	Let $G_i:=G|_{((Y_i))}$. Then the groups $G_i$ are closed, and $\partition_i:=(\emptyset, X_i^2, \smalleq|_{X_i})$ is an $\omega$-partition of $G_i$ with every component isomorphic to $N_i$. By Corollary~\ref{wreath_product} it follows that $G_i$ contains a group $G_i^*$ such that $\partition_i$ is also an $\omega$-partition of $G_i^*$ and $G_i^*$ is isomorphic to $N_i\wr \sym(\omega)$. By renaming the elements of $X$ we can assume that $G_i^*=N_i\wr \sym(\omega)$, and $K=Y_0\times \{0\}$ for some set $Y_0$ which is disjoint from $\bigcup_{i=1}^kY_i$. Then $G^*:=\id(Y_0)\times \prod_{i=1}^k(N_i\wr \sym(\omega))=\smallgroup(\id(Y_0),N_1,\dots,N_k)$ is a subgroup of $G$, and $\partition$ is an $\omega$-partition of $G^*$ with the same component as $G$. This finishes the proof of the lemma.
\end{proof}

\begin{theorem}\label{thm:main}
	$\gone_n=\gone_n'$.
\end{theorem}

\begin{proof}
	Lemma~\ref{easy_containment} shows the containment $\gone_n'\subset \gone_n$.

	Now we show $\gone_n\subset \gone_n'$ by induction on $n$. For $n=-1$ there is nothing to prove. For the induction step let us assume that $\gone_{n-1}=\gone_{n-1}'$, and let $G\in \gone$ with an $\omega$-partition $\partition=(K, \bigeq, \smalleq)$. Let $G^*$ be a subgroup as in Lemma~\ref{wreath_product}. We can assume without loss of generality that $G^*=\smallgroup(\id(Y_0),N_1,\dots,N_k)$ where $N_i\in \gone_{n-1}=\gone_{n-1}'$. It follows from the construction of $G^*$ that $G_{((Y_i\times \{0\}))}=\{(a,0)\mapsto (\sigma(a),0)\colon\sigma\in N_i\}$ where $Y_i=\dom(N_i)$. Let us observe that the equivalence relations $\bigeq$ and $\smalleq$ as defined in Definition~\ref{def:groups2} coincide with the relations $\bigeq$ and $\smalleq$ coming from the partition $\partition$. In particular $G$ preserves these equivalence relations, and also preserves the set $K:=Y_0\times \{0\}$. Therefore we can apply Lemma~\ref{recursive_description} for the groups $G^*=\smallgroup(\id(Y_0),N_1,\dots,N_k)$ and $G$. We obtain that $G=\group(H;N_0,\dots,N_k)$ for some group $H$ with 
$\prod_{i=0}^kN_i\triangleleft H$. Therefore in order to prove that $G\in \gone_n'$ it remains to show that $N,H\in \gone_{n-1}'=\gone_{n-1}$. By Lemma~\ref{direct_gone} we know that $N:=\prod_{i=1}^kN_i\in \gone_{n-1}$. Since $\gone_{n-1}$ is closed under taking closed supergroups it follows that $\overline{H}\in \gone_{n-1}$. The same argument as before shows that $N$ is also a normal subgroup of $\overline{H}$. Then by Corollary~\ref{finite_index} it follows that $|\overline{H}:N|<\infty$. Hence $|H:N|$ is also finite, and thus in fact $H$ is closed, and thus $H=\overline{H}\in \gone_{n-1}$.
\end{proof}

	We give yet another description of the class $\gone_n$ which we will need later in Section~\ref{sect:fin_bounded}.

\begin{theorem}\label{desc2}
	Let $n\in \omega$. Then a permutation group $G$ is contained in $\gone_n$ if and only if there exists a sequence of groups $G_0,G_1,\dots,G_m$ such that
\begin{enumerate}
\item $G_0=\id(\{\emptyset\})$,
\item for all $0\leq i\leq m$ one of the following holds
\begin{enumerate}
\item $G_i\simeq G_j$ for some $j<i$,
\item $G_i$ is a direct product of some groups $G_j\colon j<i$,
\item\label{item:wr} $G_i=G_j\wr \sym(\omega)$ for some $j<i$,
\item $G_i$ is a finite index supergroup of some $G_j\colon j<i$
\end{enumerate}
\item for every sequence of indices $i_0<\ldots<i_d$ for which $G_{i_l}=G_{i_{l-1}}\wr \sym(\omega)$ holds we have $d\leq n$.
\end{enumerate}
\end{theorem}

	Theorem~\ref{desc2} has the following straightforward consequence.

\begin{corollary}\label{desc2_gen}
	A permutation group $G$ is contained in $\gone$ if and only if there exists a sequence of groups $G_0,G_1,\dots,G_m$ which satisfy Conditions (1) and (2) in Theorem~\ref{desc2}.
\end{corollary}

	Corollary~\ref{desc2_gen} means in other words that $\gone$ is the smallest class of permutations groups which contains $G_0=\id(\emptyset)$ and closed under isomorphisms, finite direct products, wreath products with $\sym(\omega)$, and finite index supergroups.

	Using Theorems~\ref{reducts_groups} and~\ref{lachlan} Corollary~\ref{desc2_gen} can also be written in terms of structures as follows.

\begin{theorem}\label{desc_struct}
	$\ms$ is the smallest class of structures which contains a one-element structure, is closed under isomorphisms, and closed under taking finite disjoint unions, infinitely many copies and finite index reducts.  
\end{theorem}

\begin{remark}
	We already know that the class $\ms$ is closed under taking reducts. Therefore by taking all reducts instead of finite index reducts in the formulation of Theorem~\ref{desc_struct} we still obtain the same class of structures.
\end{remark}

\begin{proof}[Proof of Theorem~\ref{desc2}]
	First we show the ``if'' direction by induction on $m$. For $m=0$ we have $G_0=\id(\{\emptyset\})\in \gone_0$ so in this case we are done. Let us assume now that $G_0,G_1,\dots,G_m$ satisfy the conditions (1)-(3).

	\emph{Case A. $G_m=G_j$ for some $j<m$.} By the induction hypothesis $G_j\in \gone_n$. Then we also have $G_m\in \gone_n$ since $\gone_n$ is closed under isomorphisms.

	\emph{Case B. $G_m$ is a finite direct product of some groups $G_j\colon j<m$.} In this case it follows from the induction hypothesis that $G_j\in \gone_n$. Then Lemma~\ref{direct_gone} implies that $G_m\in \gone_n$. 

	\emph{Case C. $G_m=G_j\wr \sym(\omega)$ for some $j<m$.} In this case we can apply the induction hypothesis for $j<m$ and $n-1$. We obtain that $G_j\in \gone_{n-1}$. Then Corollary~\ref{inf_copies2} implies that $G_m\in \gone_n$.

	\emph{Case D. $G_m$ is a finite index supergroup of some $G_j\colon j<m$.} It follows again by the induction hypothesis that $G_j\in \gone_n$. Then $G_j$ is closed, and thus so is $G_m$. We know that the class $\gone_n$ is closed under taking closed supergroups. Therefore $G_m\in \gone_n$.

	Now we show the ``only if'' direction by induction on $n$. 

	If $n=0$ then $G$ is finite. Then we can assume that $G\leq \sym(k)$ for some $k\in \omega$. Then the sequence of groups $G_i:=\id(\{i\})\colon i<k$, $G_k=\id(k)=\prod_{i=0}^{k-1}(\id(\{i\})$, $G_{k+1}=G$ satisfy conditions (1)-(3).

	For the induction step let us assume that $G=\group(H;N_0,\dots,N_k)$ for some $H,N_0,\dots,N_k,N\in \gone_{n-1}$. Then by the induction hypothesis we know that there is a sequence $G_0,\dots,G_m$ which satisfies conditions (1),(2) and (3) for $n-1$ and which contains all groups $H,N_0,\dots,N_k$. Then let 
\begin{align*}G_{m+1}:=&N_i, G_{m+i+1}:=N_i\wr \sym(\omega),\\
 G_{m+k+1}:=&\smallgroup(N_0,\dots,N_k), G_{m+k+2}:=G=\group(H;N_1,\dots,N_k).
\end{align*}

	Then clearly the groups $G_0,\dots,G_{m+k+1}$ satisfy conditions (1),(2) and condition (3) with $n$. The only thing we have to check that we can add $G=G_{m+k+2}$ to this sequence. For this it is enough to show that $|G:\smallgroup(N_0,\dots,N_k)|$ is finite. Indeed by Lemma~\ref{small_normal_subgroup} we know that the subgroup $\smallgroup(N_0,\dots,N_k)$ is normal, and thus by Lemma~\ref{finite_index} its index is finite in $G$.
\end{proof}

\subsection{The rank 1 case}\label{sect:rank1}

	Next we focus our attention on the class $\ms_1$. Using Corollary~\ref{lachlan2} it follows that $\fa\in \ms_1$ if and only if $\aut(\fa)\in \gone_1$. We call such structures \emph{cellular}. Using Theorem~\ref{thm:main} and Remark~\ref{rank0} we know that $G\in \gone_1$ if and only if $G$ can be written, up to isomorphism, as $G=\group(H;N_0,\dots,N_k)$ for some groups $N_0,\dots,N_k,H\supset N:=\bigcup_{i=0}^kN_i$ with finite domain.

	In this section we collect some equivalent characterizations of a structure being cellular. We first show a few characterizations of cellular structures in terms of orbit growth on subsets.

\begin{definition}
	We say that a group $G$ is \emph{$P$-oligomorphic} if $\os_n(\fa)<cn^k$ for some $c>0$ and $k\in \omega$.
\end{definition}

	The class of closed $P$-oligomorphic groups are completely classified in~\cite{falque2019classification}. The following two lemmas follow from the results of~\cite{falque2019classification} and~\cite{braunfeld2019monadic}.

\begin{lemma}\label{poligomorph}
	Let $G$ be a permutation groups acting on a countable set. Then the following are equivalent.
\begin{enumerate}
\item\label{it:poli1} $G$ is $P$-oligomorphic.
\item\label{it:poli2} $\os_n(\fa)\sim cn^k$ for some $c>0$ and $k\in \omega$.
\item\label{it:poli3} $\os_n(\fa)<e^{c\sqrt{n}}$ for all $c>0$.
\end{enumerate}
\end{lemma}

\begin{lemma}\label{cell_poli}
	A structure $\fa$ is cellular if and only if it is stable and its automorphism group is $P$-oligomorphic.
\end{lemma}

	We can also characterize cellular structure in terms of orbit growth on injective tuples as follows.

\begin{lemma}\label{kexpp}
	A structure $\fa$ is cellular if and only if $\oi_n(\fa)<cn^{dn}$ for some constants $c,d$ with $d<1$.
\end{lemma}

	The advantage of the characterization above that here we do not have to assume stability, it is automatically implied by the orbit growth condition.
	
	In order to prove Lemma~\ref{kexpp} we first describe cellularity of a structure in terms of \emph{finite coverings}. Then Lemma~\ref{kexpp} is a direct consequence of the results of~\cite{bodirsky2018structures}.

\begin{definition}\label{finite_def}
	Let $\fa$ and $\fb$ be structures. A mapping $\pi\colon \fa\rightarrow \fb$ is called a \emph{finite cover}
	 if
\begin{enumerate}
\item $\pi$ is surjective,
\item for each $w\in B$ the set $\pi^{-1}(w)$ is finite,
\item the equivalence relation $\{(a,b)\colon \pi(a)=\pi(b)\}$ is preserved by $\aut(\fa)$, 
\item the image of $\aut(\fa)$ under the induced homomorphism $\iota_\pi\colon\aut(\fa)\rightarrow \sym(B)$ is equal to $\aut(\fb)$.
\end{enumerate}
	A structure $\fa$ is called a \emph{finite covering of $\fb$} if there is a finite cover $\pi\colon\fa\rightarrow \fb$.
\end{definition}

\begin{definition}
	A structure is called \emph{unary} if all of its relations are unary.
\end{definition}

	Using these definition we know the following.

\begin{theorem}\label{main:kexpp}
	Let $\fa$ be a structure. Then the following are equivalent.
\begin{enumerate}
\item\label{it:fru} $\fa$ is a finite covering of a reduct of some unary $\omega$-categorical structure.
\item\label{it:rfu} $\fa$ is a reduct of a finite covering of some unary $\omega$-categorical structure.
\item\label{it:kexpp} $\oi_n(\fa)<cn^{dn}$ for some constants $c,d$ with $d<1$.
\end{enumerate}
\end{theorem}

\begin{proof}
	See~\cite{bodirsky2018structures}, Theorem 8.16, item 8.7.
\end{proof}

	Now we show that a structure $\fa$ satisfies condition (2) of Theorem~\ref{main:kexpp} if and only if it is cellular. Together with Theorem~\ref{main:kexpp} this implies Lemma~\ref{kexpp}.

\begin{lemma}\label{m1_fcover}
	A structure $\fa$ is cellular is and only if it is a reduct of finite covering of some $\omega$-categorical unary structure.
\end{lemma}

\begin{proof}
	First, let us assume that $\fa$ is cellular. Then we can assume without loss of generality that the automorphism group of $\fa$ can be written as $\aut(\fa)=\group(H;N_0,\dots,N_k)$ for some finite domain permutation groups $N_0,\dots,N_k$ and $H\supset N:=\bigcup_{i=0}^kN_i$. Let $\fa'$ be a structure with $\aut(\fa')=\smallgroup(N_0,\dots,N_k)$. Then it is enough to show that $\fa'$ is a finite covering of some unary $\omega$-categorical structure. Let $\pi$ the map which maps every element to its $\smalleq$ class. We identify the image of $\pi$ with $\{(0,0)\}\cup \{1,\dots,k\}\times \omega$ as in Section~\ref{sect:build}. Then the image of $\aut(\fa')$ under the homomorphism $\iota(\pi)\colon\aut(\fa')\rightarrow \sym(\{(0,0)\}\cup \{1,\dots,k\}\times \omega)$ induced by $\pi$ is $\id(\{(0,0)\})\times \sym(\omega)^{[k]}$. This is exactly the automorphism group of the unary structure $\fb$ with domain set $A/\smalleq$ and with relations for the subsets $\{(0,0)\}$ and $\{i\}\times \omega$ for $i=1,\dots,k$. Clearly $\fb$ is $\omega$-categorical which finishes the proof of the ``only if'' direction of the lemma.

	Now let us assume that $\fa$ is a reduct of a finite covering $\fa'$ of some unary $\omega$-categorical structure $\fb$. We claim that $\fa$ is cellular. We have seen that the class of cellular structures is closed under taking reducts, therefore it is enough to show that $\fa'$ is cellular. We know that $\aut(\fb)$ can be written as $\prod_{i=1}^n\sym(U_i)$ for some finite partition $U_1,\dots,U_n$ of $B$ (see for instance Lemma 3.1 in~\cite{bodirsky2018structures}). We can assume without loss of generality that $U_1,\dots,U_k$ are infinite, and $U_{k+1},\dots,U_n$ are finite. Let $\pi\colon\fa'\rightarrow \fb$ be a finite covering map. We define
\begin{itemize}
\item $K:=\pi^{-1}(\bigcup_{i=k+1}^nU_i)$,
\item $\bigeq:=\{(u,v)\colon \exists i\leq k(\pi(u),\pi(v))\in U_i\}$, and
\item $\smalleq:=\{(u,v)\colon \pi(u)=\pi(v)\in \bigcup_{i=1}^kU_i\}$.
\end{itemize}
	Then it is straightforward to verify that $\partition:=(K,\bigeq,\smalleq)$ is an $\omega$-partition of $\aut(\fa')$ all of whose components having finite domains. By definition this means exactly that $\fa'$ is cellular.
\end{proof}

	Finally, we mention that in~\cite{braunfeld2019mutual} the authors showed that a structure is cellular if and only if it is $\omega$-categorical, and \emph{mutually algebraic}. We skip the definition of mutual algebraicity here, and refer the reader to the aformentioned paper.

	We summarize the results mentioned in this section in the following theorem.

\begin{theorem}\label{ms_orbit_growth}
	Let $\fa$ be a countable structure. Then the following are equivalent.
\begin{enumerate}
\item $\fa$ is cellular.
\item $\fa$ is a finite covering of reduct of a unary $\omega$-categorical structure.
\item $\fa$ is a reduct of a finite covering of a unary $\omega$-categorical structure.
\item $\aut(\fa)\in \gone_1$.
\item $\aut(\fa)\simeq \group(H;N_0,\dots,N_k)$ for some finite domain permutation groups $N_i$ and $H\supset \prod_{i=0}^kN_i$.
\item $\fa$ is stable, and $\os_n(\fa)<cn^k$ for some $c>0$ and $k\in \omega$.
\item $\fa$ is stable, and $\os_n(\fa)\sim cn^k$ for some $c>0$ and $k\in \omega$.
\item $\fa$ is stable, and $\os_n(\fa)<e^{c\sqrt{n}}$ for all $c>0$.
\item $\oi_n(\fa)<cn^{dn}$ for some constants $c,d$ with $d<1$.
\item $\fa$ is $\omega$-categorical and mutually algebraic.
\end{enumerate}
\end{theorem}

\section{Homogeneity and finite boundedness}\label{sect:fin_bounded}

	In~\cite{lachlan1992} Lachlan showed that every $\omega$-categorical monadically stable is finitely homogenizable. In this article we show that in fact $\ms\subset \fbh$. Also our proof uses some more general arguments which are not specific to the class $\ms$. The proof works as follows. By Corollary~\ref{fbh_union_closed} we know that the classes $\fhom$ and $\fbh$ are closed under taking finite disjoint unions and infinite copies. We show in Section~\ref{sect:fin_index} that these classes are also closed under taking finite index reducts. Then we can apply Theorem~\ref{desc_struct} to conclude that $\ms\subset \fbh$.

\subsection{Finite index reducts of homogeneous structures}\label{sect:fin_index}

\begin{lemma}\label{same_big_orbit}
	Let $\fa$ be an $\omega$-categorical structure, and let $\fb$ be a reduct of $\fa$. Let $\alpha_i\colon i\in I$ be such that $\aut(\fb):=\bigcup_{i\in I}(\aut(\fa)\alpha_i)$. Then $\vec{u}\same_{\fb}\vec{v}$ iff $\vec{u}\same_{\fa}\alpha_i(\vec{v})$ for some $i\in I$.
\end{lemma}

\begin{proof}
	Since $\fa$ and thus also $\fb$ are $\omega$-categorical we have
\begin{multline*}
	\vec{u}\same_{\fb}\vec{v} \Leftrightarrow \exists (\beta\in \aut(\fb))(\vec{u}=\beta(\vec{v})) \Leftrightarrow (\exists i\in I)(\exists \alpha\in \aut(\fa))(\vec{u}=\alpha\alpha_i(\vec{v}))\Leftrightarrow \\ \Leftrightarrow (\exists i\in I)(\vec{u}\same_{\fa}\alpha_i(\vec{v})).
\end{multline*}
\end{proof}

\begin{lemma}\label{finite_index_hom}
	Let $\fa$ be a homogeneous relational structure with $m(\fa)=m$, and let $\fb$ be a reduct of $\fa$ with $|\aut(\fb):\aut(\fa)|=d<\infty$. Then $\fb$ is interdefinable with a homogeneous structure $\fb'$ with $m(\fb')=dm$.
\end{lemma}

\begin{proof}
	We check condition (4) of Lemma~\ref{when_finhom} for the structure $\fb$ and $dm$. Let us assume that $\vec{u}\notsame_{\fb}\vec{v}$. We have to show that $\vec{u}_{\pi}\notsame_{\fb} \vec{v}_{\pi}$ for some $\pi\colon dm\rightarrow n$. Since $|\aut(\fb):\aut(\fa)|=d$ it follows that $\aut(\fb)$ can be written as $\aut(\fb)=\bigsqcup_{i<d}(\aut(\fa)\alpha_i)$. Then by Lemma~\ref{same_big_orbit} we know that $\vec{u}\notsame_{\fa}\alpha_i(\vec{v})$ for all $i<d$. Then by using the implication (1)$\rightarrow$(4) of Lemma~\ref{when_finhom} we obtain that for all $i<d$ there exists some $\pi_i\colon m\rightarrow n$ such that $\vec{u}_{\pi_i}\notsame_{\fa}\alpha_i(\vec{v}_{\pi_i})$. Let $$\pi\colon dm\rightarrow n, im+j\mapsto \pi_i(j), (i< d, j< m).$$ Then $\vec{u}_{\pi}\notsame_{\fa}\alpha_i(\vec{v}_{\pi})$ for all $i< d$. By using Lemma~\ref{same_big_orbit} again this implies $\vec{u}_{\pi}\notsame_{\fb} \vec{v}_{\pi}$, and this is what we wanted to show.
\end{proof}

	Lemma~\ref{finite_index_hom} immediately implies the following.

\begin{corollary}\label{finite_index_hom2}
	The class $\fhom$ is closed under taking finite index reducts.
\end{corollary}

\begin{proof}
	Let $\fa\in \fhom$, and let $\fb$ be a reduct of $\fa$ with $|\aut(\fb):\aut(\fa)|=d<\infty$. We show that $\fb\in \fhom$. We can assume without loss of generality that $\fa$ is homogeneous and finitely related. Then Lemma~\ref{finite_index_hom} immediately implies $\fb$ is interdefinable with a homogeneous finite related structure. 
\end{proof}

	Now we show that the class $\fbh$ is also closed under taking finite index reducts.

\begin{lemma}\label{how_to_merge}
	Let $\fa$ be a homogeneous finitely bounded structure with signature $\sigma$, and let $m:=b(\fa)$. Let $\tau\subset \sigma$, and let $\fb$ be a $\tau$-reduct of $\fa$ which is also homogeneous. Let $\fc$ be a $\tau$-structure, and let $a_0,\dots,a_m$ be pairwise distinct elements in $\fc$. Then $\fc\in \age(\fb)$ if and only if there exist a $\sigma$-expansion of $\fd_i$ of the structures $\fc_i:=\fc|_{C\setminus \{a_i\}}$ such that for all $i,j\leq m$ we have
\begin{itemize}
\item $\fd_i\in \age(\fa)$, and
\item $\fd_i|_{C\setminus \{a_i,a_j\}}=\fd_j|_{C\setminus \{a_i,a_j\}}$.
\end{itemize}
\end{lemma}

\begin{proof}
	We use the notation $C_i:=C\setminus \{a_i\}$, and $C_{ij}=C\setminus \{a_i,a_j\}=C_i\cap C_j$.

	Let $\fc\in \age(\fb)$. We can assume without loss of generality that $\fc$ is a finite substructure of $\fb$. Then let $\fd:=\fa|_{C}$, and $\fd_i:=\fa|_{C_i}$. Then obviously $\fd_i\in \age(\fa)$, and $\fd_i|_{C_{ij}}=\fd_j|_{C_{ij}}=\fd|_{C_{ij}}$.

	For the other direction we define a $\sigma$-structure $\fd$ as follows. The domain set of $\fd$ is $C$, and $\vec{u}\in R^{\fd}$ if and only if there exists some $i$ such that $\vec{u}\in R^{\fd_i}$.

	\emph{Claim 1. $\fd|_{C_i}=\fd_i$.} Let $R\in \sigma$. If $\vec{u}\in R^{\fd_i}$ then it follows immediately from the definition that $\vec{u}\in R^{\fd}$. On the other hand if $\vec{u}\in R^{\fd}\cap C_i^{\ar(R)}$ then by definition $\vec{u}\in R^{\fd_j}\cap C_i^{\ar(R)}$ for some $j\leq m$. Then $\vec{u}\in C_j^{\ar(R)}\cap C_i^{\ar(R)}=C_{ij}^{\ar(R)}$, and thus item (2) implies that $\vec{u}$ is also contained in the relation $R^{\fd_i}$.

	\emph{Claim 2. $\fd\in \age(\fa)$}. Let us assume that this is not the case. Then let $\fd'$ be a minimal substructure of $\fd$ which is not in age $\age(\fa)$. Then $|\fd'|\leq b(\fa)=m$. Therefore there exists some $a_i\colon 0\leq i\leq m$ which is not contained in $\fd'$. Then the domain set of $\fd'$ is contained in $C_i$ for some $i$, and thus $\fd'$ is substructure of $\fd|_{C_i}=\fd_i$. In particular $\fd'\in \age(\fa)$, a contradiction.

	Finally, let $\fc$ be the $\tau$-reduct of $\fd$. Then $\fc\in \age(\fb)$.
\end{proof}

\begin{lemma}\label{fin_types}
	Let $\fa$ be an $\omega$-categorical structure with signature $\sigma$, and let $\fb$ be a reduct of $\fa$ with signature $\tau$. Let $d:=|\aut(\fb):\aut(\fa)|$. Then for every $\tau$-type $p$ which is realizable in $\fb$ there exists at most $d$ many $\sigma$-types $q$ which are realizable in $\fa$ and such that $p\subset q$.
\end{lemma}

\begin{proof}
	Let $\aut(\fb):=\bigcup_{i\in I}(\aut(\fa)\alpha_i)$ with $|I|=d$, and let us choose $\vec{u}\in A^k$ with $\tp_{\fb}(\vec{u})=p$. Then if $q\supset p$ then there exists some $\vec{v}\in A^k$ such that $\tp_{\fa}(\vec{v})=q$ and $\tp_{\fb}(\vec{v})=p=\tp_{\fb}(\vec{u})$. By Lemma~\ref{same_big_orbit} the second conditions holds if and only if $\vec{v}\same_{\fa}\alpha_i(\vec{u})$ for some $i\in I$. This implies that if $\tp_{\fb}(\vec{v})=p$ then $q=\tp_{\fa}(\vec{v})=\tp_{\fa}(\sigma_i(\vec{u}))$ for some $i\in I$. This shows the statement of the lemma. 
\end{proof}

\begin{lemma}\label{finite_index_bounded}
	Let $\fa$ be a homogeneous finitely bounded structure with signature $\sigma$, and let $m:=b(\fa)$. Let $\tau\subset \sigma$, and let $\fb$ be a $\tau$-reduct of $\fa$ which is also homogeneous, and let us assume that $d:=|\aut(\fb):\aut(\fa)|$ is finite. Then $b(\fb)\leq md^{m+1}$.
\end{lemma}

\begin{proof}
	Let $\fc\in \Min(\comp{\age(\fb)})$. We have to show that $|\fc|\leq md^{m+1}$. Let $a_0,\dots,a_m\in C$, and let $\fc_i:=\fc|_{C\setminus \{a_i\}}$. Then $\fc_i\in \age(\fb)$. Therefore there exist $\sigma$-expansions $\fd_i$ of $\fc_i$ which are contained in $\age(\fa)$. Let $\vec{u}_i$ be a tuple which enumerates all elements of $C_i$. Let $\fd_{i0},\dots,\fd_{i(l(i)-1)}$ be all the possible $\sigma$-expansions of $\fc_i$ which are contained in $\age(\fa)$. Then $\tp_{\fd_{ij}}(\vec{u})\colon j<l(i)$ are pairwise different $\sigma$-types all of which contain $\tp_{\fc}(\vec{u})$. Since $\fd_{ij}\in \age(\fa)$ it follows that all $\sigma$-types $\tp_{\fd_{ij}}(\vec{u})$ are realizable in $\fa$. Therefore by Lemma~\ref{fin_types} it follows that $l(i)< d$ for all $i\leq m$.

	Since $\fc\not\in \age(\fb)$, and $\fd_{ij}\colon i\leq m, j<l(i)$ are all in $\age(\fa)$ it follows by Lemma~\ref{how_to_merge} that for all functions $f\colon m+1\rightarrow d$ with $f<l$ we have $\fd_{if(i)}|_{C\setminus \{a_i,a_j\}}\neq \fd_{jf(j)}|_{C\setminus \{a_i,a_j\}}$ for some $i,j\leq m$. This implies that for every function $f$ as above there exists a tuple $\vec{u}_f$ from $C_{ij}:=C_i\cap C_j$ and a $R_f\in \sigma$
such that $\vec{u}_f\in R_f^{\fd_{if(i)}}\not\Leftrightarrow \vec{u}_f\in R_f^{\fd_{jf(j)}}$. Now let $C'$ be the union of element in all tuples $\vec{u}_f$. Then length of the tuples $\vec{u}_f$ is at most $m(\fa)\leq b(\fa)\leq m$. Therefore $|C'|\leq md^{m+1}$. Let $\fc':=\fc|_{C'}$.

	Now we show that $\fc'\not\in \age(\fb)$ $(*)$. Let us assume for contradiction that $\fc'\in \age(\fb)$. Let $C_i':=C_i\setminus \{a_i\}$, $C_{ij}':=C_i'\cap C_j'$, and $\fc_i':=\fc|_{C' \setminus \{a_i\}}$. Then applying Lemma~\ref{how_to_merge} again it follows that there exists $\sigma$-expansions $\fd'_i$ of $\fc_i'$ such that $\fd_i'\in \age(\fa)$, and $\fd_i'|_{C_{ij}'}=\fd_j'|_{C_{ij}'}$. It follows from the homogeneity of $\fb$ that for all $i$ there a $\sigma$-expansion $\fd_i''\in \age(\fa)$ of $\fc_i$ such that $\fd_i''|_{C_i'}=\fd_i'$. Then $\fd_i''=\fd_{if(i)}$ for some $f(i)< l(i)$ for all $i$. This implies that for all $i,j< m$ we also have $\fd_{if(i)}|_{C_i'\cap C_j'}=\fd_{jf(j)}|_{C_i'\cap C_j'} (\dag)$. On the other by definition we know that there exist indices $i,j\leq m$, a tuple $\vec{u}_f$ from $C_{ij}$, and a relational symbol $R_f$ such that $\vec{u}_f\in R_f^{\fd_{if(i)}}\not\Leftrightarrow \vec{u}_f\in R_f^{\fd_{jf(j)}}$. By the construction of $C'$ it follows that every coordinate of $\vec{u}_f$ is from $C'$ and then also from $C'\cap C_{ij}=C'\setminus \{a_1,a_2\}=C_{ij}'$. Combining this with the observation that $\fd_{if(i)}|_{C_i'}=\fd_i'$ it follows that $R^{\fd_i'}(\vec{u}_f)\not\Leftrightarrow R^{\fd_j'}(\vec{u}_f)$. This contradicts the equality $(\dag)$, and thus finishes the proof of $(*)$.

	Since $\fc$ was chosen to be a minimal $\tau$-structure in $\comp{\age(\fb)}$, the statement $(*)$ implies that in fact $\fc=\fc'$. Hence $|\fc|=|\fc'|\leq md^{m+1}$ which finishes the proof of the lemma.
\end{proof}

\begin{lemma}\label{finite_index_bounded2}
	Let $\fa$ be a homogeneous finite bounded structure with $m:=b(\fa)$, and let $\fb$ be a reduct of $\fa$ with $|\aut(\fb):\aut(\fa)|=d$. Then $\fb$ is interdefinable with a homogeneous structure $\fb'$ with $b(\fb')\leq md^{dm+2}$.
\end{lemma}

\begin{proof}
	Let $\fb':=\typestruct_{dm}(\fb)$. Then $\fb'$ is homogeneous and it is interdefinable with $\fb$ by Lemmas~\ref{finite_index_hom} and~\ref{when_finhom}. Let $\fa'$ be the structure $\fa$ expanded by all relations of $\fb'$. Then $\fa$ and $\fa'$ are interdefinable, and thus $\aut(\fa)=\aut(\fa')$ and $\fa'$ is also homogeneous. Then by applying Lemma~\ref{finite_index_bounded} for the structures $\fa'$ and $\fb'$ we obtain $b(\fb')\leq (dm)d^{dm+1}=md^{dm+2}$.
\end{proof}

	Lemma~\ref{finite_index_bounded2} immediately implies the following.

\begin{corollary}\label{thm:finite_index}
	The class $\fbh$ is closed under taking finite index reducts.
\end{corollary}

\begin{theorem}\label{thm:ms_fbh}
	$\ms\subset \fbh$.
\end{theorem}

\begin{proof}
	Follows immediately from Theorem~\ref{desc_struct}, Corollary~\ref{fbh_union_closed}, and Corollary~\ref{thm:finite_index}.
\end{proof}

\section{Thomas' conjecture for hereditarily cellular structures}\label{sect:thomas}

	In this section we show that Thomas' conjecture holds for the class $\ms$ i.e., every structure in $\ms$ has finitely many reduct up to interdefinability.

	Let $G\leq \sym(X)$ be a permutation group. Then we define $\width(G)$ to be the supremum of the set $$\{|\acl_{G_x/E}(\emptyset)|\colon x\in X, E \text{ is a congruence of } G_x\}.$$

\begin{lemma}\label{w_finite}
	If $G$ is oligomorphic then $\width(G)$ is finite.
\end{lemma}

\begin{proof}
	If $x$ and $y$ are in the same orbit of $G$ then 
\begin{align*}
&\sup\{|\acl_{G_x/E}(\emptyset)|\colon E \text{ is a congruence of } G_x\}=\\&\sup\{|\acl_{G_y/E}(\emptyset)|\colon E \text{ is a congruence of } G_y\}.
\end{align*}
	Therefore it is enough to show that $\sup\{|\acl_{G_x/E}(\emptyset)|\colon E \text{ is a congruence of } G_x\}$ is finite for a fixed $x$. Since $G$ is oligomorphic, so is $G_x$, and thus $G_x$ has finitely many congruences. Therefore it is enough to show that $\acl_{G_x/E}(\emptyset)$. Since $G_x$ is oligomorphic, so is $G_x/E$. In particular $G_x/E$ has finitely many finite orbits, and therefore $\acl_{G_x/E}(\emptyset)$ is finite.
\end{proof}

\begin{lemma}\label{w_reduct}
	If $G\leq H$ then $\width(G)\geq \width(H)$.
\end{lemma}

\begin{proof}
	If $E$ is a congruence of $H_x$ for some $x\in X$ then it is also a congruence of $G_x$, and $G_x/E\leq H_x/E$. Therefore in this case $\acl_{H_x/E}(\emptyset)\leq \acl_{G_x/E}(\emptyset)$. Taking the supremum over $x$ and the congruences $E$ of $H_x$ the claim of the lemma follows.
\end{proof}

\begin{lemma}\label{w_product}
	Let $G_1,\dots,G_k$ be permutation groups with pairwise disjoint domains, and let $G:=\prod_{i=1}^kG_i$. Then $\width(G)\leq \sum_{i=1}^k\width(G_i)$.
\end{lemma}

\begin{proof}
	Let $X_i:=\dom(G_i)$, and $X:=\dom(G)$. Let $x\in X_i$, and let $E$ be a congruence of $G$. Let $\{[y_1]_E,\dots,[y_m]_E\}$ be a finite orbit of $G_x/E$ for some $y_1,\dots,y_m$ with $y_1\in X_i$. Then since $G_x$ fixes the sets $X_1,\dots,X_k$ it follows that we can find $y_1',\dots,y_m'\in X_j$ such that $[y_l']_E=[y_l]_E$ for all $l=1,\dots,m$. Then $\{[y_1']_{E_j},\dots,[y_m']_{E_j}\}$ is a finite orbit of $(G_x|_{X_i})/E_j$ for $E_i:=E\cap X_i^2$. This implies that $\acl_{G_x/E}(\emptyset)$ can be written as a union of $\acl_{(G_x)|_{X_j}/E_j}(\emptyset)\colon j=1\dots,k$. Since $G:=\prod_{i=1}^kG_i$ it follows that $G_x|_{X_j}=G_j$ if $j\neq i$, and $G_x|_{X_j}=G_i|_x$ if $j=1$. Let $x_j\in X_j$ be arbitrary with $x_i=x$. Then since $(G_j)_{x_j}\leq (G_x)|_{X_j}$ for all $j$, and thus $\acl_{(G_x)|_{X_j}/E_j}(\emptyset)\subset \acl_{(G_j)_{x_j}/E_j}(\emptyset)$, and thus $|\acl_{(G_x)|_{X_j}/E_j}(\emptyset)|\leq \width(G_j)$. Therefore $|\acl_{G_x/E}(\emptyset)|\leq \sum_{i=1}^k\width(G_i)$. Since this holds for every $x\in X$, and a congruence $E$ of $G_x$ we obtain $\width(G)\leq \sum_{i=1}^k\width(G_i)$.
\end{proof}

\begin{lemma}\label{counting_width}
	Let $Y_0,\dots,Y_k,Y,N_0,\dots,N_k,N,H$ be as in Section~\ref{sect:build}. Let $H_i:=H_{(Y_i)}$. Then $\width(H_i)\leq \width(\group(H;N_0,\dots,N_k))$ for all $i=0,1,\dots,k$.
\end{lemma}

\begin{proof}
	In this proof we use the notations from Section~\ref{sect:build}.
	
	Let $G:=\group(H;N_0,\dots,N_k)$. We have to show that for all $x\in Y$ and a congruence $E$ of $(H_i)_x$ we have $|\acl_{(H_i)_x/E}(\emptyset)|\leq \width(G)$.

	Let us assume that $x\in Y_0$. Then $\acl_G(\emptyset)=Y_0$, and thus $\width(G)\geq |Y_0|$. On the other hand $\width((H_0)_x)\leq |\dom(H_0)|=|Y_0|$. Therefore the statement of the lemma if true for $i=0$. 

	Let $x\in Y_i$ for $i>0$, and let $E$ be a congruence of $(H_i)_x$. Let $E^*$ the equivalence relation on $X$ whose equivalence classes are

\begin{enumerate}[(a)]
\item $Y_0\times \{0\}$,
\item all $\bigeq$-classes except $Y_i\times \omega$,
\item all $\smalleq$-classes in $Y_i$ except $Y_i\times \{0\}$,
\item the set of the form $C\times \{0\}$ where $C\in Y_i/E$. 
\end{enumerate}	

	We claim that $E^*$ is a congruence of $G_{(x,0)}$. Since $G_{(x,0)}$ preserves the relations $\bigeq$ and $\smalleq$ it follows that $G$ maps every $E^*$ class corresponding to cases (a)-(d) to another class corresponding to the same item. Let $C\in Y_i/E_x$, and let $\sigma\in G_{(x,0)}$ we have to show that $\sigma$ maps $C\times \{0\}$ to $C'\times \{0\}$ for some $C'\in Y_i/E_x$. Let $\vec{0}:=(0,\dots,0)$. Then by definition $\rho_{\vec{0}}(\sigma)\in H$. Since $G_{(x,0)}$ preserves $(x_0)$ and the relation $\smalleq$ it follows that $\sigma$ preserves $[(x,0)]_{\smalleq}=Y_i\times 0$. Therefore $\sigma(a,0):=(\rho_{\vec{0}}(\sigma),0)$ for all $a\in Y_i$. This also implies that $\rho_{\vec{0}}(\sigma)\in Y_i$, and thus $\sigma':=\rho_{\vec{0}}(\sigma)|_{Y_i}\in H_{(Y_i)}=H_i$. Since $\sigma((x,0))=(x,0)$ we have $\sigma'(x)=x$. Then $C':=\sigma'(C)\in Y_i/E$ since $E$ is a congruence of $(H_i)_x$, and thus $$\sigma(C\times \{0\})=((\rho_{\vec{0}}(\sigma)(C))\times \{0\}=\sigma'(C)\times \{0\}=C'\times \{0\}.$$ This shows that $E^*$ is a congruence of $E$.

	The definition of $E^*$ implies that $\acl_{(H_i)_x/E}(\emptyset)\times \{0\}\subset \acl_{G_{(x,0)}/E^*}(\emptyset)$, and thus $|\acl_{(H_i)_x/E}(\emptyset)|\leq |\acl_{G_{(x,0)}/E^*}(\emptyset)|\leq \width(G)$ which finishes the proof of the lemma.
\end{proof}

\begin{lemma}\label{subdirect_direct}
	Let $N_0,\dots,N_k,N,H\in \gone$ be as in Definition \ref{def:groups2}, and let $H_i:=H_{(Y_i)}$. Then the group $H^*:=\langle H, \prod_{i=0}^kH_i\rangle$ is closed, and $N$ is a finite index normal subgroup of $H^*$.
\end{lemma}

\begin{proof}
	By our assumption $N$ is a normal subgroup of $H$, therefore it is enough to show that $N$ is normalized by $\prod_{i=0}^kH_i$. Since $N\triangleleft H$ it follows that $\id(Y_0)\times \dots \id(Y_{i-1})\times N_i \times \id(Y_{i+1})\dots \times \id(Y_k)$ is normalized by $H|_{\{Y_i\}}$. Therefore $N_i$ is normalized by $H_i$, and thus $N=\prod_{i=1}^kN_i$ is normalized by $\prod_{i=0}^kH_i$. This also means that $N\triangleleft \overline{H^*}$. By Corollary~\ref{finite_index} this implies that $|H^*:N|\leq |\overline{H^*}:N|<\infty$, and thus $H^*$ is closed.
\end{proof}

\begin{lemma}\label{counting_width2}
	Let $n,m\in \omega$. Then there exists finitely many closed groups $G\in \gone_n$ with $\width(G)\leq m$ up to isomorphism.
\end{lemma}

\begin{proof}
	We prove the lemma by induction on $n$. For $n=-1$ there is nothing to show. 
	
	For the induction step let us consider a group $G\in \gone_n$ with $\width(G)\leq m$ for some $m\in \omega$. By Theorem~\ref{thm:main} we know that $G$ is isomorphic to group $\group(H;N_0,\dots,N_k)$ where $N_0,N_1,\dots,N_k,H\in \gone_{n-1}$. Let $Y_i:=\dom(N_i)$, $H_i:=H_{(Y_i)}$, $H^*:=\langle H, \prod_{i=0}^kH_i\rangle$, and $\smalleq^*:=\bigcup_{i=0}^kY_i^2$. Then $\smalleq^*$ is a congruence of $G$, and $\dom(G/\smalleq^*)\geq k+1$. Therefore $k+1\leq \width(G)=m$. By Lemma~\ref{subdirect_direct} we know that $N\triangleleft H^*$ and $H^*$ is closed. Hence $H^*\in \gone_{n-1}$. By Lemma~\ref{counting_width} it follows that $\width(H_i)\leq m$, and then Lemma~\ref{w_product} implies $\width(\prod_{i=0}^kH_i)\leq m(k+1)$. By Lemma~\ref{w_reduct} we have $\width(H^*)\leq \width(\prod_{i=0}^kH_i)$, therefore $\width(H^*)\leq m(k+1)\leq m^2$. By the induction hypothesis this implies that we have finitely many choices for the group $H^*$ up to isomorphism. This means that it is enough to show that for a fixed group $H^*\in \gone_{n-1}$ we have finitely many choices for groups $N$ and $H$ in $\gone_{n-1}$ such that $N\triangleleft H^*$ and $H\leq H^*$. By Corollary~\ref{finite_index} we know that we have finitely many possible choices for the group $N$, and in all cases $|H^*:N|$ is finite. Therefore if both $H^*$ and $N$ are fixed then we have finitely many choices for $N\leq H\leq H^*$. This finishes the proof of the lemma.
\end{proof}

\begin{corollary}\label{countable}
	$\gone$ contains countable many groups up to isomorphism.
\end{corollary}

\begin{proof}
	Direct consequence of Lemma~\ref{counting_width2}.
\end{proof}

	Now we can show that Thomas' conjecture holds for the class $\ms$.

\begin{theorem}\label{thomas}
	Every hereditarily cellular structure has finitely many reduct up to interdefinability.
\end{theorem}

\begin{proof}
	Let $\fa\in \ms$. We will check Condition (5) of Proposition~\ref{thomas_equivalent} $\fa$. Let $G:=\aut(\fa)$. Then $G\in \gone$ by Theorem~\ref{lachlan}, and thus $G\in \gone_n$ for some $n\in \omega$. Let $m:=\width(G)$. Then for every closed supergroup $H$ we have $H\in \gone_n$ and $\width(H)\leq m$ by Lemma~\ref{w_reduct}. Therefore by Lemma~\ref{counting_width2} we have finitely many choices for $H$ up to isomorphism.
\end{proof}

\begin{remark}
	We know by Theorem~\ref{thm:ms_fbh} (or by the results of~\cite{lachlan1992}) that structure in $\ms$ is interdefinable with a finitely related homogeneous structure. This implies that Theorem~\ref{thomas} is indeed a special case of Thomas' conjecture.
\end{remark}

\bibliographystyle{alpha}
\bibliography{local.bib}

\end{document}